\documentclass[a4paper]{article}
\usepackage{RR}
\usepackage{hyperref}

\usepackage[english]{babel}
\usepackage[latin1]{inputenc}
\usepackage[T1]{fontenc}
\usepackage{float}
\usepackage{amsmath}
\usepackage{amsfonts}
\usepackage{amssymb}
\usepackage{amsthm}
\usepackage{stmaryrd}
\usepackage{url}
\usepackage{epsfig}
\usepackage{natbib}
\usepackage{color}

\def\var{\mbox{var}}
\def\corr{\mbox{corr}}
\def\AIC{\mbox{AIC}}
\def\BIC{\mbox{BIC}}
\def\pen{\mbox{pen}}

\newenvironment{itemize*}%
  {\begin{itemize}%
    \setlength{\itemsep}{0pt}%
    \setlength{\parskip}{0pt}}%
  {\end{itemize}}

\psfigdriver{dvips}
\usepackage{mathrsfs}

\usepackage{geometry}
\newfloat{Figure}{H}{lof}
\newfloat{Table}{Hptb}{lot} 

\newtheorem{thrm}{Theorem}[section]
\newtheorem{prte}[thrm]{Proposition}
\newtheorem{lemma}[thrm]{Lemma}

\newtheorem{defi}{Definition}[section]

\newtheorem{algo}{Algorithm}[section]

\RRdate{Janvier 2009}
\RRversion{2}
\RRdater{Septembre  2009}

\RRauthor{Nicolas Verzelen
\thanks{Laboratoire de Math\'ematiques UMR 8628, Universit\'e Paris-Sud, 91405 Osay}
\thanks{INRIA Saclay, Projet SELECT, Universit\'e Paris-Sud, 91405 Osay}}
\authorhead{Verzelen}

\RRtitle{Sélection automatique de voisinage d'un champ gaussien}
\RRetitle{Data-driven neighborhood selection of a Gaussian field}
\titlehead{Neighborhood selection}
\RRresume{Nous étudions l'estimation non-paramétrique d'un champ gaussien stationnaire $X$ observé sur un réseau régulier. Dans ce cadre, nous avons précédemment introduit une procédure de sélection de modèle \cite{verzelen_gmrf_theorie}. Cette procédure revient à sélectionner un voisinage $\widehat{m}$ grâce une technique de pénalisation puis à estimer la covariance du champ $X$ dans l'espace des champs de Markov gaussiens de voisinage $\widehat{m}$. Une telle stratégie satisfait des inégalités oracles et des propriétés d'apdaptation au sens minimax. En pratique, elle présente néanmoins quelques inconvénients. D'une part, la pénalité dépend de quantités inconnues. D'autre part, la procédure est uniquement définie pour des réseaux toriques. La contribution de cet article est triple. Nous proposons un algorithme automatique pour calibrer la pénalité. De plus, nous introduisons une extension à des réseaux non-toriques. Enfin, nous étudions les performances pratiques de la procédure sur des données simulées. Ces simulations suggèrent que la sélection de champs de Markov gaussiens est souvent une bonne alternative à l'estimation de variogramme.}

\RRabstract{We study the nonparametric covariance estimation of a stationary Gaussian field $X$ observed on a lattice. To tackle this issue, a neighborhood selection procedure has been recently introduced. This procedure amounts to selecting a neighborhood $\widehat{m}$ by a penalization method and estimating the covariance of $X$ in the space of Gaussian Markov random fields (GMRFs) with neighborhood $\widehat{m}$. Such a strategy is shown to satisfy oracle inequalities as well as minimax adaptive properties. However, it suffers several drawbacks which make the method difficult to apply in practice. On the one hand, the penalty depends on some unknown quantities. On the other hand, the procedure is only defined for toroidal lattices.
The present contribution is threefold. A data-driven algorithm is proposed for tuning the penalty function. Moreover, the procedure is extended to non-toroidal lattices. Finally, numerical study illustrate the performances of the method on simulated examples. These simulations suggest that Gaussian Markov random field selection is often a good alternative to variogram estimation.}
\RRmotcle{Champ gaussien, champ de Markov gaussien, calibration automatique, sélection de modèle, pseudo-vraisemblance.}
\RRkeyword{Gaussian field, Gaussian Markov random field, Data-driven calibration, model selection, pseudolikelihood.}

\RRprojets{Select}
\RRtheme{\THCog} 
 \RCSaclay 

\begin{document}

\RRNo{6798}
\makeRR   

\section{Introduction}\label{introduction}


We study the  estimation of the distribution of a stationary Gaussian field $(X{\scriptstyle [i,j]})_{(i,j)\in\Lambda}$ indexed by the nodes of a rectangular lattice $\Lambda$ of size $p_1\times p_2$.  This problem is often encountered in spatial statistics or in image analysis. Classical statistical procedures allow 
to estimate and subtract the trend. Henceforth, we assume that the field $X$ is centered. Given a $n$-sample of the field $X$, the challenge is to  infer the correlation. In practice, the number $n$ of observations often equals one. Different methods have been proposed to tackle this problem. 

A traditional approach amounts to computing an empirical variogram and then fitting a suitable parametric variogram model such as the exponential or Mat\'ern model (see \cite{cressie} Ch.2 or \cite{stein}). The main disadvantage with this method is that the practitioner is required to select  a \emph{good} variogram model. When the field exhibits long range dependence, specific procedures have been introduced (e.g. Fr\'ias\emph{ et al.}~\cite{frias07}). In the sequel, we focus on small range dependences. Most of the nonparametric (Hall \emph{et al.}~\cite{hall94}) and semiparametric (Im\emph{ et al.}~\cite{stein07}) methods are based on the spectral representation of the field. To our knowledge, these procedures have not yet been shown to achieve adaptiveness, i.e. their rate of convergence does not adapt to the \emph{complexity} of the correlation functions. 

In this paper, we define and study a nonparametric estimation procedure relying on Gaussian Markov random fields (GMRF). This procedure is computationally fast and satisfies adaptive properties. 
Let us fix a node  $(0,0)$ at the center of $\Lambda$   and let $m$ be a subset of $\Lambda\setminus\{(0,0)\}$. The field $X$ is a GMRF with respect to the neighborhood $m$ if conditionally to $(X{\scriptstyle[k,l]})_{(k,l)\in m}$, the variable $X{\scriptstyle[0,0]}$ is independent from all the remaining variables in $\Lambda$. We refer to Rue and Held \cite{rue} for a comprehensive introduction on GMRFs. If we know that $X$ is a GMRF with respect to the neighborhood $m$, then we can estimate the covariance by applying  
likelihood or pseudolikelihood maximization. Such parametric procedures are  well understood, at least from an asymptotic point of view  (see for instance \cite{guyon95} Sect.4). However,  we do not know in practice what is the ``good'' neighborhood $m$.
For instance, choosing the empty neighborhood amounts to assuming that all the components of $X$ are independent. Alternatively, if we choose  the complete neighborhood, which contains all the nodes of $\Lambda$ except $(0,0)$, then the number of parameters is huge and estimation performances are poor.\\

We tackle in this paper the problem of neighborhood selection from a practical point of view. The purpose is to define a data-driven procedure that picks a suitable neighborhood $\widehat{m}$ and then estimates the distribution of $X$ in the space of GMRFs with neighborhood $\widehat{m}$. This procedure neither requires	 any knowledge on the correlation of $X$, nor assumes that the field $X$ satisfies a Markov condition. Indeed, the procedure selects a neighborhood $\widehat{m}$ that achieves a trade-off between an \emph{approximation} error (distance between the true correlation and GMRFs with neighborhood $m$) and an \emph{estimation} error (variance of the estimator). If  $X$ is a GMRF with respect to a small neighborhood, then the procedure achieves a parametric rate of convergence. Alternatively, if $X$ is not a GMRF then the rate of convergence of the procedure depends on the rate of approximation of the true covariance by GMRFs with growing neighborhood. In short, the procedure is nonparametric and adaptive.

Besag and Kooperberg \cite{besag95}, Rue and Tjelmeland \cite{rue02}, Song \emph{et al.} \cite{fuentes2008}, and Cressie and Verzelen \cite{verzelen_cressie} have considered  the problem of \emph{approximating} the correlation of a Gaussian field by a GMRF, but this approach requires the knowledge of the true distribution. Guyon and Yao have stated in \cite{guyon2000} necessary conditions and sufficient conditions for a model selection procedure to choose asymptotically the true neighborhood of a GMRF with probability one. Our point of view is slightly different.
We do not assume that the field $X$ is a GMRF with respect to a sparse neighborhood. We do not aim at estimating the true neighborhood, we rather want to select a neighborhood that allows to estimate \emph{well} the distribution of $X$ (i.e. to minimize a risk). The distinction between these two points of view has been nicely described in the first chapter of MacQuarrie and Tsai \cite{MacQuarrie}.\\

In \cite{verzelen_gmrf_theorie}, we have introduced a neighborhood selection procedure based on pseudolikelihood maximization and penalization. Under mild assumptions, the procedure achieves optimal neighborhood selection. More precisely, it satisfies an oracle inequality and it is  minimax adaptive to the sparsity of the neighborhood. To our knowledge, these are the first results of neighborhood selection in this spatial setting.

If the procedure exhibits appealing theoretical properties, it  suffers several drawbacks from a practical perspective.
First, the method constrains the largest eigenvalue of the estimated covariance to be smaller than some parameter $\rho$. In practice, it is difficult to choose  $\rho$ since we do not know the largest eigenvalue of the true covariance. Second, the penalty function $\pen(.)$ introduced in Sect.3 of the previous paper  depends on the largest eigenvalue of the covariance of the field $X$.
Hence, we need a practical method for tuning the penalty. Third,  the procedure has only been defined when the lattice $\Lambda$ is a square torus. \\

	
Our contribution is twofold. On the one hand,  we propose practical versions of our neighborhood selection procedure that overcome the previously-mentioned drawbacks:
\begin{itemize*}
 \item The procedure is extended to rectangular lattices.
\item We do not constrain anymore the largest eigenvalue of the covariance.
\item  We provide an algorithm based on the  so-called \emph{slope heuristics} of Birg\'e and Massart \cite{massart_pente} for tuning the penalty. Theoretical justifications for its use are also given.
\item Finally, we extend the procedure to the case where the lattice $\Lambda$ is not a torus. 
\end{itemize*}
On the other hand, we illustrate the performances of this new procedure on numerical examples.
When $\Lambda$ is a torus, we compare it with likelihood-based methods like $\AIC$ \cite{akaike73} and $\BIC$ \cite{BIC}, even if they were not studied in this setting. When $\Lambda$ is not toroidal, likelihood methods become intractable. Nevertheless, our procedure still applies and often outperforms variogram-based methods. \\

The paper is organized as follows. In Section \ref{section_presentation}, we define a new version of the estimation procedure of \cite{verzelen_gmrf_theorie} that does not require anymore  the choice of the constant $\rho$.  We also discuss the computational complexity of the procedure. In Section \ref{section_resultats_theorique}, we connect this new procedure to the original method and we recall some theoretical results. We provide an algorithm for tuning the penalty in practice in Section  \ref{section_slope_heuristic}. In Section \ref{section_extension_non_torique}, we extend our procedure for handling non-toroidal lattices. The simulation studies are provided in Section \ref{section_simulations}. Section \ref{section_conclusion} summarizes our findings, while  the proofs are postponed to Section \ref{section_proofs}.\\

Let us introduce some notations. In the sequel, $X^v$ refers  to the vectorialized version of $X$ with the convention $X{\scriptstyle[i,j]}=X^v{\scriptstyle[(i-1)\times p_2+ j]}$ for any $1\leq i\leq p_1$ and $1\leq j\leq p_2$. Using this new notation amounts to ``forgetting'' the spatial structure of $X$ and allows to get into a more classical statistical framework. We note ${\bf X}_1, {\bf X}_2,\ldots ,{\bf X}_n$ the $n$ observations of the field $X$. The matrix $\Sigma$ stands for the covariance matrix of $X^v$.  For any matrix $A$,  $\varphi_{\text{max}}(A)$ and  $\varphi_{\text{min}}(A)$ respectively refer the largest eigenvalue and the smallest eigenvalues of $A$. Finally, $I_{r}$ denotes the identity matrix of size $r$.

\section{Neighborhood selection on a torus}\label{section_presentation}

In this section, we introduce the main concepts and notations for GMRFs on a torus. Afterwards, we describe our procedure based on pseudolikelihood maximization. Finally, we discuss some computational aspects. Throughout this section and the two following sections, the lattice $\Lambda$ is assumed to be toroidal. Consequently, the components of the matrices $X$  are taken modulo $p_1$ and $p_2$.

\subsection{GMRFs on the torus}

The notion of conditional distribution is underlying the definition of GMRFs. By standard Gaussian derivations (see for instance \cite{lauritzen} App.C), there exists a unique $p_1\times p_2$ matrix $\theta$ such that $\theta{\scriptstyle[0,0]}=0$ and
\begin{eqnarray}\label{regression_conditionnelle}	
X{\scriptstyle[0,0]}= \sum_{(i,j)\in \Lambda\setminus\{(0,0)\}}\theta{\scriptstyle[i,j]}X{\scriptstyle[i,j]} + \epsilon{\scriptstyle[0,0]}\ ,
\end{eqnarray}
where the random variable $\epsilon{\scriptstyle[0,0]}$ follows a zero-mean normal distribution and is independent from the covariates $(X{\scriptstyle[i,j]})_{(i,j)\in \Lambda\setminus\{(0,0)\}}$. The linear combination $\sum_{(i,j)\in\Lambda\setminus\{(0,0)\}}\theta{\scriptstyle[i,j]}X{\scriptstyle[i,j]}$ is the kriging predictor of $X{\scriptstyle[0,0]}$ given the remaining variables. In the sequel, we note $\sigma^2$ the variance of $\epsilon{\scriptstyle[0,0]}$ and we call it the conditional variance of $X{\scriptstyle[0,0]}$.

Equation (\ref{regression_conditionnelle}) describes the conditional distribution of $X{\scriptstyle[0,0]}$ given the remaining variables. By stationarity of the field $X$, it holds that that $\theta{\scriptstyle[i,j]}=\theta{\scriptstyle[-i,-j]}$. The covariance matrix $\Sigma$ is closely related to $\theta$ through the following equation:  
\begin{eqnarray}\label{lien_sigma_theta}
 \Sigma = \sigma^2 \left[I_{p_1p_2} - C(\theta)\right]^{-1} \ ,
\end{eqnarray}
where the $p_1p_2\times p_1p_2$ matrix $C(\theta)$ is defined by $C(\theta){\scriptstyle[(i_1-1)p_2+j_1,(i_2-1)p_2+j_2]} := \theta{\scriptstyle[i_2-i_1,j_2-j_1]}$ for any $1\leq i_1,i_2\leq p_1$ and $1\leq j_1,j_2\leq p_2$. The matrix $(I_{p_1p_2} - C(\theta))$ is called the partial correlation matrix of the field $X$.
The so-defined matrix $C(\theta)$ is symmetric block circulant with $p_2\times p_2$ blocks. We refer to \cite{rue} Sect.2.6 or the book of Gray \cite{gray} for definitions and main properties on circulant and block circulant matrices.\\

Identities (\ref{regression_conditionnelle}) and (\ref{lien_sigma_theta}) have two main consequences. First, estimating the $p_1\times p_2$ matrix $\theta$ amounts to estimating the covariance matrix $\Sigma$ up to a multiplicative constant. We shall therefore focus on $\theta$. Second, by Equation (\ref{regression_conditionnelle}), the field $X$ is a GMRF with respect to the neighborhood defined by the support $\theta$.
The adaptive estimation issue of the distribution of $X$ by neighborhood selection therefore reformulates as an adaptive estimation problem of the matrix $\theta$ via support selection.

Let us now precise the set of possible values for $\theta$. The set $\Theta$ denotes the vector space of the $p_1\times p_2$ matrices that satisfy $\theta{\scriptstyle[0,0]}=0$ and 
$\theta{\scriptstyle[i,j]}=\theta{\scriptstyle[-i,-j]}$, for any $(i,j)\in\Lambda$. Hence, a matrix $\theta\in\Theta$ corresponds to the distribution of a stationary Gaussian field if and only if the $p_1p_2\times p_1p_2$ matrix $(I_{p_1p_2}-C(\theta))$ is positive definite. This is why we define the convex subset $\Theta^+$ of $\Theta$ by
\begin{eqnarray}
 \Theta^+:=\left\{\theta\in\Theta\hspace{0.2cm}\text{s.t.}\,  \left[I_{p_1p_2}-C(\theta)\right]\text{ is positive definite}\right\}\label{definition_theta+}\ .
\end{eqnarray}
The set of covariance matrices of stationary Gaussian fields on $\Lambda$ with unit conditional variance is in one to one correspondence with the set $\Theta^+$. We sometimes assume that the field $X$ is isotropic. 
The corresponding sets $\Theta^{\text{iso}}$ and $\Theta^{+,\text{iso}}$ for isotropic fields are introduced as: 
\begin{eqnarray*}
 \Theta^{\text{iso}} :=\left\{\theta\in\Theta\ ,\hspace{0.2cm}\theta{\scriptstyle[i,j]}=\theta{\scriptstyle[-i,j]}=\theta{\scriptstyle[j,i]}\ ,\,  \forall (i,j)\in\Lambda \right\}\ \text{ and } \Theta^{+,\text{iso}} := \Theta^+\cap\Theta^{\text{iso}}\ .
\end{eqnarray*}
	
\subsection{Description of the procedure}
Let $|(i,j)|_t$ refer to the toroidal norm defined by 
\begin{eqnarray*}
|(i,j)|^2_t :=\left[i\wedge \left(p_1-i\right)\right]^2 +\left[j\wedge \left(p_2-j\right)\right]^2\ ,
\end{eqnarray*}
for any node $(i,j)\in\Lambda$.

In the sequel, a model $m$ stands for a subset of $\Lambda\setminus\{(0,0)\}$. It is  also called a neighborhood.
For the sake of simplicity, we shall only use the  collection of models $\mathcal{M}_1$ defined below. 
\begin{defi}\label{definition_modele}
A subset $m\subset \Lambda\setminus\{(0,0)\}$ belongs to $\mathcal{M}_{1}$ if and only if 
 there exists a number $r_m>1$ such that 	
\begin{eqnarray}
m = \left\{ (i,j)\in \Lambda\setminus\{(0,0)\}\hspace{0.3cm} \text{s.t.}\hspace{0.3cm}  |(i,j)|_t\leq r_m\right\}\ . \label{definition_rm}
\end{eqnarray}
\end{defi}
In other words, the neighborhoods $m$ in $\mathcal{M}_1$ are sets of nodes lying in a disc centered at $(0,0)$.
Obviously, $\mathcal{M}_{1}$ is totally ordered with respect to the
inclusion. Consequently, we order the models $m_0\subset m_1\subset \ldots \subset m_{i}\ldots  $. For instance, $m_0$ corresponds to the empty neighborhood, $m_1$ stands for the neighborhood of size $4$, and $m_2$ refers to the neighborhood with $8$ neighbours. See Figure \ref{figure_modele} for an illustration.

\begin{Figure}
\centerline{
{\bf a)} \epsfig{file=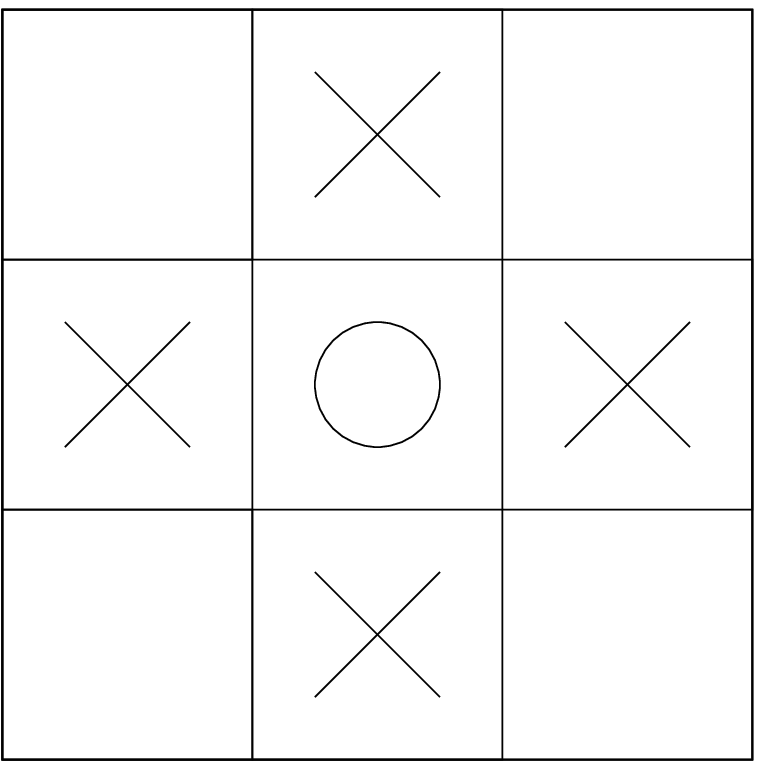,width=3cm}\hspace{0.6cm}
{\bf b)} \epsfig{file=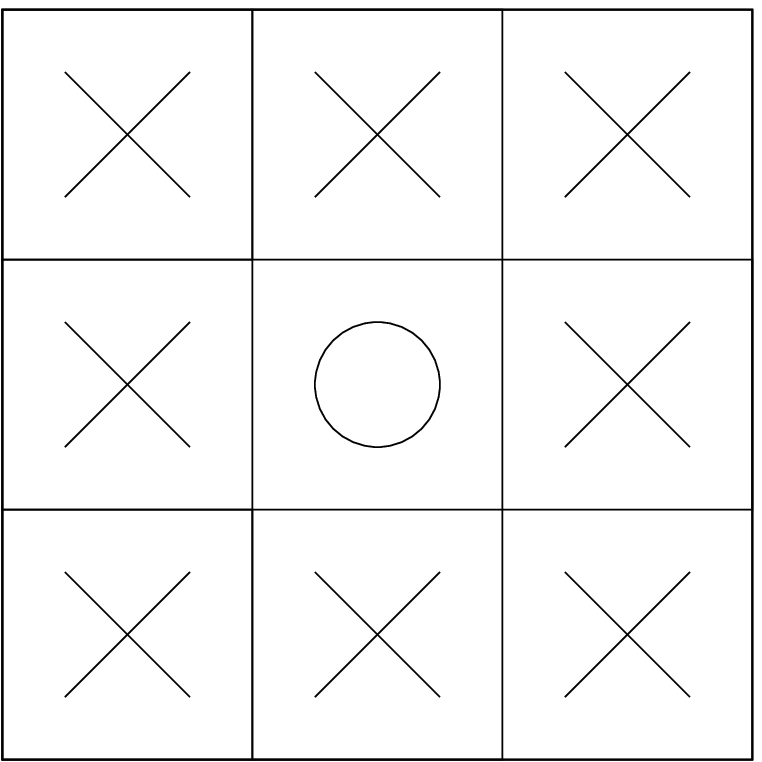,width=3cm}\hspace{0.6cm}
{\bf c)} \epsfig{file=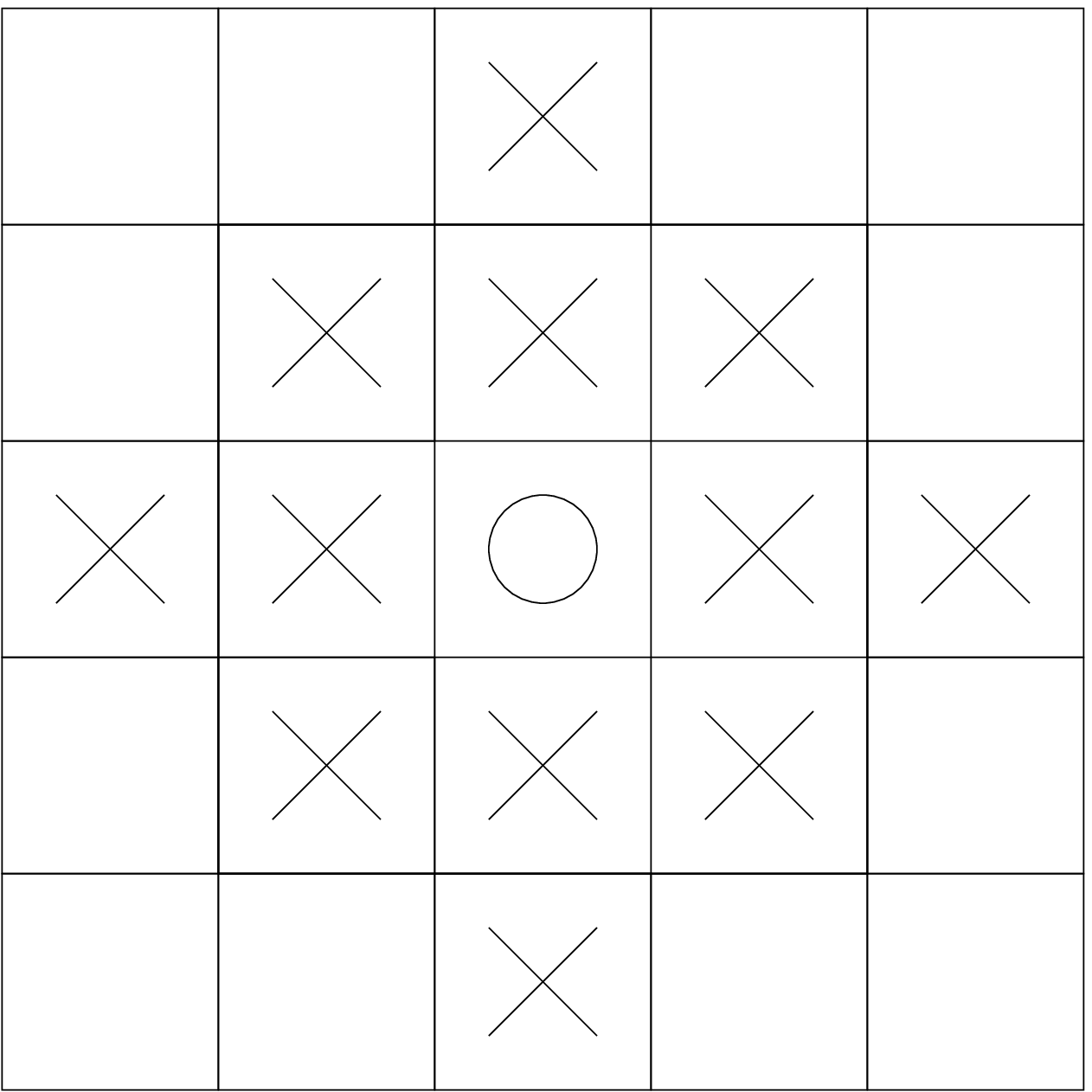,width=5cm}}
\caption{\textit{(a) Model $m_1$ with first order neighbors. (b) Model $m_2$ with second order neighbors. (c) Model $m_3$ with third order neighbors. }}
\label{figure_modele}
\end{Figure}

For any model $m\in\mathcal{M}_1$, the vector space $\Theta_m$ is  the subset of matrices $\Theta$ whose support is included in $m$. Similarly $\Theta^{\text{iso}}_m$ is the subset of $\Theta^{\text{iso}}$ whose support is included in $m$. The dimensions of $\Theta_m$  and $\Theta^{\text{iso}}_m$ are respectively noted $d_m$ and $d_m^{\text{iso}}$. 
Since we aim at estimating the positive matrix $(I_{p_1p_2}-C(\theta))$, we also consider the convex subsets
 of $\Theta_m^+$ and $\Theta_m^{+,\text{iso}}$ which correspond to non-negative precision matrices.
\begin{eqnarray}\label{definition_thetam+}
 \Theta_m^+ :=\Theta_m\cap \Theta^+ \hspace{1.5cm}\text{ and }\hspace{1.5cm} \Theta_m^{+,\text{iso}} :=\Theta^{\text{iso}}_m\cap \Theta^{+,\text{iso}}\ .
\end{eqnarray}

For any $\theta'\in \Theta^+$, the conditional least-squares (CLS)  criterion  $\gamma_{n,p_1,p_2}(\theta')$~	\cite{guyon}  is defined by
\begin{eqnarray}\label{definition_gamman}
\gamma_{n,p_1,p_2}(\theta') & := & \frac{1}{np_1p_2}\sum_{i=1}^n\sum_{(j_1,j_2)\in\Lambda} \bigg({\bf X}_i{\scriptstyle[j_1,j_2]} -
\sum_{(l_1,l_2)\in \Lambda\setminus\{(0,0)\}} \theta'{\scriptstyle[l_1,l_2]}{\bf X}_i{\scriptstyle[j_1+l_1,j_2+l_2]}\bigg)^2 \ .
\end{eqnarray}
The function $\gamma_{n,p_1,p_2}(.)$ is a least-squares criterion that allows us to perform the simultaneous linear regression of all ${\bf X}_i{\scriptstyle[j_1,j_2]}$ with respect to the covariates $({\bf X}_i{\scriptstyle[l_1,l_2]})_{(l_1,l_2)\neq(k_1,k_2)}$.  This criterion is closely connected with the pseudolikelihood introduced by Besag \cite{besag75-2}. The associated estimator is slightly less efficient estimator than maximum likelihood estimation (\cite{guyon95} Sect.4.3). Nevertheless, its computation is much faster since it does not involve determinants as for the likelihood. See \cite{verzelen_gmrf_theorie} Sect. 7.1, for a more complete comparison between CLS and maximum likelihood estimators in this setting. For any model $m\in\mathcal{M}_1$, the estimators are defined as the unique minimizers of $\gamma_{n,p_1,p_2}(.)$ on the sets $\Theta_{m}^+$ and $\Theta_m^{+,\text{iso}}$.
\begin{eqnarray}
 \widehat{\theta}_{m}:=\arg\min_{\theta'\in\overline{\Theta_{m}^+}} \gamma_{n,p_1,p_2}(\theta')\hspace{1.5cm}\text{and}\hspace{1.5cm}  \label{definition_estimateur} \widehat{\theta}^\text{iso}_{m}:=\arg\min_{\theta'\in\overline{\Theta_{m}^{+,\text{iso}}}} \gamma_{n,p_1,p_2}(\theta')\ , 
\end{eqnarray}
where $\overline{A}$ stands for the closure of $A$.  We further discuss the connection between $\widehat{\theta}_{m}$ and $\widehat{\theta}_{m,\rho_1}$ in Section \ref{section_resultats_theorique}.

Given a subcollection of models $\mathcal{M}$ of $\mathcal{M}_1$ and a positive function $\pen:\mathcal{M}\rightarrow \mathbb{R}^+$ called a penalty, we select a model as follows:
\begin{eqnarray}\label{definition_modele_selection}
 \widehat{m} := \arg\min_{m\in\mathcal{M}}\left[\gamma_{n,p_1,p_2}\left(\widehat{\theta}_m\right)+\pen(m)\right]\hspace{0.3cm}\text{and}\hspace{0.3cm}\widehat{m}^{\text{iso}} := \arg\min_{m\in\mathcal{M}}\left[\gamma_{n,p_1,p_2}\left(\widehat{\theta}^{\text{iso}}_m\right)+\pen(m)\right]\ .
\end{eqnarray}
For short, we write $\widetilde{\theta}$ and $\widetilde{\theta}^{\text{iso}}$ for $\widehat{\theta}_{\widehat{m}}$ and $\widehat{\theta}_{\widehat{m}^{\text{iso}}}^{\text{iso}}$. We discuss the choice of the penalty function  in Section \ref{section_slope_heuristic}.

\subsection{Computational aspects}\label{section_computationelle}

Since the lattice $\Lambda$ is a torus, the computation of the estimators $\widehat{\theta}_m$ is performed efficiently thanks to the following lemma.

\begin{lemma}\label{lemme_transformee_fourier}
For any $p\times p$ matrix $A$ and for any $1\leq i\leq p_1$  and $1\leq j\leq p_2$, let $\lambda{\scriptstyle[i,j]}(A)$ be the $(i,j)$-th term of two-dimensional discrete Fourier transform of the matrix $A$, i.e. 
\begin{eqnarray}\label{definition_fft}
 \lambda{\scriptstyle[i,j]}(A) &:= &\sum_{k=1}^{p_1}\sum_{l=1}^{p_2}A{\scriptstyle[i,j]} \exp\bigg[2\iota\pi\left(\frac{ki}{p_1}+\frac{jl}{p_2}\right)\bigg] \ ,
\end{eqnarray}
where $\iota^2=-1$. The conditional least-squares criterion $\gamma_{n,p_1,p_2}(\theta')$ simplifies as  
	\begin{eqnarray*} 
\gamma_{n,p_1,p_2}(\theta') & = &\frac{1}{np_1^2p_2^2}\bigg\{\sum_{i=1}^{p_1}\sum_{j=1}^{p_2}\left[1-\lambda{\scriptstyle[i,j]}(\theta)\right]^2\bigg[\sum_{k=1}^n\lambda{\scriptstyle[i,j]}\left({\bf X}_k\right)\overline{\lambda{\scriptstyle[i,j]}\left({\bf X}_k\right)}\bigg]\bigg\}\ .
\end{eqnarray*}
\end{lemma}
A proof is given in Section \ref{section_proofs}. Optimization of $\gamma_{n,p_1,p_2}(.)$ over the set $\Theta_{m}^+$ is performed fastly using the fast Fourier transform (FFT).  Nevertheless, this is not the privilege of  CLS estimators, since maximum likelihood estimators are also computed fastly by FFT when $\Lambda$ is a torus.

In Section \ref{section_extension_non_torique}, we mention that the computation of the CLS estimators $\widehat{\theta}_m$ remains quite easy when $\Lambda$ is not a torus whereas likelihood maximization becomes intractable.

\section{Theoretical results}\label{section_resultats_theorique}
 Throughout this section,  $\Lambda$ is assumed to be a toroidal square lattice and we note $p$ its size. Let us mention that the restriction to square lattices made in \cite{verzelen_gmrf_theorie} allows to simplify the proofs but is not necessary so that the theoretical results hold.
In this section, we first recall the original procedure and we emphasize the differences with the one defined in the previous section. We also mention a result of optimality. This will provide  some insights for calibrating the penalty $\pen(.)$ in Section \ref{section_slope_heuristic}.\\

Given $\rho>2$ be a positive constant, we define the subsets $\Theta_{m,\rho}^{+}$ and $\Theta_{m,\rho}^{+,\text{iso}}$ by
\begin{eqnarray}
\Theta_{m,\rho}^{+} & := &\left\{\theta\in  \Theta_m^{+}\ ,\ \varphi_{\text{max}}\left[I_{p_1p_2}-C(\theta)\right]< \rho\right\} \label{defintion_thetam+gamma}\\
\Theta_{m,\rho}^{+,\text{iso}}  & := &\left\{\theta\in  \Theta_m^{+,\text{iso}}\ ,\	 \varphi_{\text{max}}\left[I_{p_1p_2}-C(\theta)\right]< \rho\right\} \ . \nonumber
\end{eqnarray}
Then, the corresponding estimators $\widehat{\theta}_{m,\rho}$ and $\widehat{\theta}^{\text{iso}}_{m,\rho}$ are defined as in (\ref{definition_estimateur}), except that we now consider $\Theta_{m,\rho}^+$ instead of $\Theta_{m}^+$.
Let us mention that the estimator $\widehat{\theta}_m$ corresponds to 
the estimator $\widehat{\theta}_{m,\rho_1}$  defined in \cite{verzelen_gmrf_theorie} Sect.2.2 with $\rho_1=+\infty$.

\begin{eqnarray*}
 \widehat{\theta}_{m,\rho}:=\arg\min_{\theta'\in\overline{\Theta_{m,\rho}^+}} \gamma_{n,p,p}(\theta')\hspace{1.5cm}\text{and}\hspace{1.5cm}  \widehat{\theta}^\text{iso}_{m,\rho}:=\arg\min_{\theta'\in\overline{\Theta_{m,\rho}^{+,\text{iso}}}} \gamma_{n,p,p}(\theta')\ .
\end{eqnarray*}
Given a subcollection $\mathcal{M}$ of $\mathcal{M}_1$  and a penalty function $\pen(.)$, we select the  models $\widehat{m}_{\rho}$ and $\widehat{m}^{\text{iso}}_{\rho}$ as in (\ref{definition_modele_selection}) except that we use $\widehat{\theta}_{m,\rho}$ and $\widehat{\theta}^{\text{iso}}_{m,\rho}$ instead of $\widehat{\theta}_{m}$ and $\widehat{\theta}^{\text{iso}}_{m}$. We also note $\widetilde{\theta}_{\rho}$ and $\widetilde{\theta}^{\text{iso}}_{\rho}$ for $\widehat{\theta}_{\widehat{m}_{\rho},\rho}$ and $\widehat{\theta}^{\text{iso}}_{\widehat{m}^{\text{iso}}_{\rho},\rho}$.\\

The only difference between the estimators $\widetilde{\theta}$ and $\widetilde{\theta}_{\rho}$ is that the largest eigenvalue of the precision matrix $(I_{p^2}-C(\widetilde{\theta}))$ is restricted to be smaller than $\rho$. We make this restriction in \cite{verzelen_gmrf_theorie} to facilitate the analysis. \\

In order to assess the performance of the penalized estimator $\widetilde{\theta}_{\rho}$ and $\widetilde{\theta}^{\text{iso}}_{\rho}$, we use the prediction loss function $l(\theta_1,\theta_2)$ defined by
\begin{eqnarray}\label{definition_perte}
l(\theta_1,\theta_2) := \frac{1}{p^2}tr\left[(C(\theta_1)-C(\theta_2))\Sigma(C(\theta_1)-C(\theta_2))\right]\ .
\end{eqnarray}
As explained in \cite{verzelen_gmrf_theorie} Sect.1.3, the loss $l(\theta_1,\theta_2)$ expresses in terms of conditional expectation
\begin{eqnarray}\label{definition_perte_alternative}
 l(\theta_1,\theta_2) =  \mathbb{E}_{\theta}\left\{\left[\mathbb{E}_{\theta_1}\left(X{\scriptstyle[0,0]}|X_{\Lambda\setminus\{0,0\}}\right)-\mathbb{E}_{\theta_2}\left(X{\scriptstyle[0,0]}|X_{\Lambda\setminus\{0,0\}}\right)\right]^2\right\}\ ,
\end{eqnarray}
where $\mathbb{E}_{\theta}(.)$ stands for the expectation with respect to the distribution $\mathcal{N}(0,\sigma^2(I_{p_1p_2}-C(\theta))^{-1})$. Hence, $l(\widehat{\theta},\theta)$ corresponds the mean squared prediction loss of $X{\scriptstyle [0,0]}$ given the other covariates. A similar loss function is also used by Song \emph{et al.} \cite{fuentes2008}, when approximation Gaussian fields by GMRFs.
For any neighborhood $m\in\mathcal{M}$, we define the \emph{projection} $\theta_{m,\rho}$ as the closest element of $\theta$ in  $\Theta_{m,\rho}^+$ with respect to the loss $l(.,.)$. 
\begin{eqnarray*}
\theta_{m,\rho}:=\arg\min_{\theta'\in\overline{\Theta_{m,\rho}^+}} l(\theta',\theta)\hspace{1.5cm}\text{and}\hspace{1.5cm}  \theta^\text{iso}_{m,\rho}:=\arg\min_{\theta'\in\overline{\Theta_{m,\rho}^{+,\text{iso}}}} l(\theta',\theta)\ .
\end{eqnarray*}
We call the loss $l(\theta_{m,\rho},\theta)$ the bias of the set $\Theta_{m,\rho}^+$. This implies  that  $\widehat{\theta}_{m,\rho}$ cannot perform better than this loss.

\begin{thrm}\label{mainthrm}
Let $\rho>2$,  $K$ be a positive number larger than an universal constant $K_0$ and $\mathcal{M}$ be a subcollection of $\mathcal{M}_1$. If for every model $m\in \mathcal{M}$, it holds that 
\begin{eqnarray} \label{condition_penalite}
\pen(m) \geq K\rho^2\varphi_{\text{max}}(\Sigma)\frac{d_m+1}{np^2}\ ,
\end{eqnarray}
then for any $\theta\in\Theta^+$ , the estimator $\widetilde{\theta}_{\rho}$  satisfies
\begin{eqnarray}\label{majoration_risque}
\mathbb{E}_{\theta}[l(\widetilde{\theta}_{\rho},\theta)] \leq
L(K)\inf_{m\in \mathcal{M}}\left[l(\theta_{m,\rho},\theta) + \pen(m)\right] \ ,
\end{eqnarray}
where $L(K)$ only depends on $K$. A similar bound holds if one replaces $\widetilde{\theta}_{\rho}$ by  $\widetilde{\theta}_{\rho}^{\text{\emph{iso}}}$, $\Theta^+$ by $\Theta^{+,\text{iso}}$, $\theta_{m,\rho}$ by $\theta_{m,\rho}^{\text{\emph{iso}}}$, and $d_m$ by $d_m^{\text{\emph{iso}}}$.
\end{thrm}
Although we have assumed the correlation is non-singular, the theorem still holds if the spatial field is constant. The nonasymptotic bound is provided in a slightly different version in \cite{verzelen_gmrf_theorie}. It states that $\widetilde{\theta}_{\rho}$ achieves a trade-off between the bias and a variance term if  the penalty is suitable chosen.  In Theorem \ref{mainthrm}, we use the penalty $K\rho^2\varphi_{\text{max}}(\Sigma)(d_m+1)/(np^2)$ instead of the penalty $K\rho^2\varphi_{\text{max}}(\Sigma)d_m/(np^2)$ stated in the previous paper. This makes the bound (\ref{majoration_risque}) simpler. Observe that these two penalties yield the same model selection since they only differ by a constant.
Let us further discuss two points.

\begin{itemize}
 \item In this paper, we use the estimator $\widetilde{\theta}$ rather than $\widetilde{\theta}_{\rho}$. Given a collection of models $\mathcal{M}$, there exists some finite $\rho>2$, such that these two estimators coincide. Take for  instance $\rho=\sup_{m\in\mathcal{M}}\sup_{\theta\in\Theta_m^+}\varphi_{\max}(I_{p_1p_2}-C(\theta))$. Admittedly, the so-obtained $\rho$ may be large, especially if there are large models in $\mathcal{M}$. The upper bound (\ref{majoration_risque}) on the risk therefore becomes worse. Nevertheless, we do not think that the dependency of (\ref{majoration_risque}) on $\rho$ is sharp. Indeed , we illustrate in Section \ref{section_simulations} that the risk of  $\widetilde{\theta}$ exhibits good statistical performances.

\item Theorem \ref{mainthrm} provides a suitable form of the penalty for obtaining oracle inequalities. 
However, this penalty depends on $\varphi_{\text{max}}(\Sigma)$ which is not known in practice. This is why we develop a data-driven penalization method in the next section.
\end{itemize}

\section{Slope Heuristics}\label{section_slope_heuristic}

Let us introduce a data-driven method for calibrating the penalty function $\pen(.)$. It is based on the so-called \emph{slope heuristic} introduced by Birg\'e and Massart \cite{massart_pente}  in the fixed design Gaussian regression framework (see also \cite{massartflour} Sect.8.5.2). This heuristic relies on the notion of minimal penalty. In short, assume that one knows that a good penalty has a form $\pen(m)=NF(d_m)$ (where $d_m$ is the dimension of the model and $N$ is a tuning parameter). Let us define $\widehat{m}(N)$ the selected model as a function of $N$. There  exists a quantity $\widehat{N}_{\text{min}}$ satisfying the following property: If $N>\widehat{N}_{\text{min}}$, the dimension of the selected model $d_{\widehat{m}(N)}$ is reasonable and if $N<\widehat{N}_{\text{min}}$, the dimension of the selected model is huge. The function $pen_{\text{min}}(.):=\widehat{N}_{\text{min}} F(.)$ is called the minimal penalty. In fact, a
\emph{dimension jump} occurs for $d_{\widehat{m}(N)}$ at the point $\widehat{N}_{\text{min}}$. Thus, the quantity $\widehat{N}_{\text{min}}$ is clearly observable for real data sets. In their Gaussian framework, Birg\'e and Massart have shown that twice the minimal penalty is nearly the  optimal penalty. In other words, the model $\widehat{m}:=\widehat{m}(2\widehat{N}_{\text{min}})$ yields an efficient estimator.

The slope heuristic method has been    successfully applied for multiple change-point detection \citep{lebarbier05}. Applications are also being developed in other frameworks such as mixture models \citep{michel07}, clustering \citep{baudry07}, estimation of oil reserves \citep{lepez02}, and genomic \citep{villers07}.

If this method was originally introduced for fixed design Gaussian regression, Arlot and Massart \citep{arlot_pente} have proved more recently that a similar phenomenon occurs in the heteroscedastic random-design case. In the GMRF setting,  we are only able to partially justify  this heuristic. For the sake of simplicity, let us assume in the next proposition that the lattice $\Lambda$ is a square of size $p$.

\begin{prte}\label{proposition_pente}
Consider $\rho>2$, and $\eta<1$ and suppose that $p$ is larger than some numerical constant $p_0$. Let $m'$ be  the largest model in $\mathcal{M}_1$ that satisfies $d_{m'}\leq \sqrt{np^2}$. 
For any model $m\in \mathcal{M}_1$, we assume that
\begin{eqnarray}
\pen(m') - \pen(m) \leq K_1(1-\eta)\sigma^2\left\{\varphi_{\text{min}}\left(I_{p^2}-C(\theta)\right)\wedge\left[\rho-\varphi_{\text{max}}\left(I_{p^2}-C(\theta)\right)\right] \right\}\frac{d_{m'}-d_m}{np^2} \ ,  \label{hypothese_pente}
\end{eqnarray}
where $K_1$ is a universal (constant defined in the proof).
Then, for any  $\theta\in\Theta^+_{m',\rho}$, it holds that
$$\mathbb{P}\left\{d_{\widehat{m}_{\rho}}> L\left[\sqrt{np^2}\wedge p^2 \right]\right\}\geq \frac{1}{2}\ ,$$
where   $L$ only depends on $\eta$, $\rho$, $\varphi_{\text{min}}\left(I_{p^2}-C(\theta)\right)$ , and $\varphi_{\text{max}}\left(I_{p^2}-C(\theta)\right)$.
\end{prte}
The proof is postponed to Section \ref{section_proofs}. Let us define
$$N_1:=K_1\sigma^2\left\{\varphi_{\text{min}}\left(I_{p_1p_2}-C(\theta)\right)\wedge\left[\rho-\varphi_{\text{max}}\left(I_{p_1p_2}-C(\theta)\right)\right]\right\} \ ,$$
and let us consider penalty functions $\pen(m)=N\frac{d_m}{np_1p_2}$ for some $N>0$.
The proposition states that if $N$ is smaller than $N_1$, then the procedure selects a model of huge dimension with large probability, i.e $d_{\widehat{m}(N)}$ is huge. Alternatively, let us define
$$N_2:=K_0\frac{\sigma^2\rho^2}{\varphi_{\text{min}}\left(I_{p_1p_2}-C(\theta)\right)}\frac{d_m}{np_1p_2}\ ,$$
where the numerical constant $K_0$ is introduced in Theorem 3.1 in  \cite{verzelen_gmrf_theorie}. By Theorem 3.1, choosing $N> N_2$ ensures that the risk of $\widetilde{\theta}_{\rho}$ achieves a type-oracle inequality and the dimension $d_{\widehat{m}_{\rho}(N)}$ is reasonable. The quantities $N_1$ and $N_2$ are different especially when the eigenvalues of $(I_{p_1p_2}-C(\theta))$ are far from $1$. Since we do not know the behavior of the selected model $\widehat{m}_{\rho}(N)$ when $N$ is  between  $N_1$ and $N_2$, we are not able to really prove a dimension jump as the fixed design Gaussian regression framework. Besides, we have mentioned in the preceding section that we are more interested in the estimator $\widetilde{\theta}$ than $\widetilde{\theta}_{\rho}$. Nevertheless,  we clearly observe in simulation studies a dimension jump for some $N$ between $N_1$ and $N_2$
even if we use the estimators $\widehat{\theta}_m$ instead of $\widehat{\theta}_{m,\rho}$.
This suggests that the slope heuristic is still valid in the GMRF framework.

\begin{algo}(Data-driven penalization with slope heuristic).\label{pente}
Let $\mathcal{M}$ be a subcollection of $\mathcal{M}_1$.
\begin{enumerate}
\item Compute the selected model $\widehat{m}(N)$ as a function of $N>0$
$$\widehat{m}(N) \in \arg\min_{m\in\mathcal{M}}\left\{\gamma_{n,p_1,p_2}\left(\widehat{\theta}_m\right)+ N \frac{d_m}{np_1p_2}\right\}\ .$$ 
\item Find $\widehat{N}_{\text{\emph{min}}}>0$ such that the jump $d_{\widehat{m}\left(\left[\widehat{N}_{\text{\emph{min}}}\right]_-\right)}-d_{\widehat{m}\left(\left[\widehat{N}_{\text{\emph{min}}}\right]_+\right)}$ is maximal. 
\item Select the model  $\widehat{m}= \widehat{m}(2\widehat{N}_{\text{\emph{min}}})$.
\end{enumerate}
\end{algo}
The difference $f(x_-)-f(x_+)$ measures the discontinuity of a function $f$ at the point $x$.
Step 2 may need to introduce huge models in the collection $\mathcal{M}$ all the other ones being considered as ``reasonably small''. 
As the function $\widehat{m}(.)$ is piecewise linear with at most $\text{Card}(\mathcal{M})$ jumps, so that steps 1-2 have a complexity $\mathcal{O}\left(\text{Card}(\mathcal{M})\right)^2$. We refer to App.A.1 of \cite{arlot_pente} for more details on the computational aspects of steps $1$ and $2$. Let us mention that there are other ways of estimating $\widehat{N}_{\text{min}}$ than choosing the largest jump as described in \cite{arlot_pente} App.A.2. 
Finally, the methodology described in this section straightforwardly extends to the case of isotropic GMRFs estimation by replacing $\widehat{m}(N)$ by $\widehat{m}^{\text{iso}}(N)$ and $d_m$ by $d_{m}^{\text{iso}}$.\\

In conclusion, the neighborhood selection procedure described in Algorithm \ref{pente} is completely data-driven and does not require any prior knowledge on the matrix $\Sigma$. Moreover, its computational burden remains small. We illustrate  its efficiency in Section \ref{section_simulations}.

\section{Extension to non-toroidal lattices}\label{section_extension_non_torique}

It is often artificial to consider the field $X$ as stationary on a torus. However, we needed this hypothesis for deriving nonasymptotic properties of the estimator $\widetilde{\theta}$ in \cite{verzelen_gmrf_theorie}. In many applications, it is more realistic to assume that we observe a small window of a Gaussian field defined on the plane $\mathbb{Z}^2$. If we are unable to prove nonasymptotic risk bounds in this new setting. Nevertheless,  Lakshman and Derin have shown in \citep{Lakshman93} that there is no phase transition within the valid parameter space for GMRFs defined on the plane $\mathbb{Z}^2$. Let us briefly explain what this means:
consider a GMRF defined on a square lattice of size $p$, but only observed on a square lattice of size $p'$. The absence of phase transition implies the distribution of this field observed on this fixed window of size $p'$ does not asymptotically depend on the bound conditions when $p$ goes  to infinity. Consequently, it is reasonable  to think that our estimation procedure still performs well 
to the price of slight modifications. In the sequel, we assume that the field $X$ is defined on $\mathbb{Z}^2$, but the data ${\bf X}$ still correspond to $n$ independent observations of the field $X$ on the window $\Lambda$ of size $p_1\times p_2$.  The conditional distribution of $X{\scriptstyle[0,0]}$ given the remaining covariates now decomposes as
\begin{eqnarray} \label{regression_plan}
 X{\scriptstyle[0,0]} = \sum_{(i,j)\in \mathbb{Z}^2\setminus\{(0,0)\}}\theta{\scriptstyle[i,j]}X{\scriptstyle[i,j]}+\epsilon{\scriptstyle[0,0]}\ ,
\end{eqnarray}
where $\theta{\scriptstyle[.,.]}$ is an ``infinite'' matrix defined on $\mathbb{Z}^2$ and where $\epsilon{\scriptstyle[0,0]}$ is a centered Gaussian variable of variance $\sigma^2$ independent of $(X{\scriptstyle[i,j]})_{(i,j)\in\Lambda\setminus\{(0,0)\}}$. The distribution of the field $X$ is uniquely defined by the function $\theta$ and positive number $\sigma^2$. The set $\Theta^{+,\infty}$ of valid parameter for $\theta$ is now defined using the spectral density function. We refer to Rue and Held \cite{rue} Sect.2.7 for more details.

\begin{defi}\label{definition_Theta_infini}
A function $\theta:\mathbb{Z}^2\rightarrow \mathbb{R}$ belongs to the set $\Theta^{+,\infty}$ if  it satisfies the three following conditions:
\begin{enumerate}
\item $\theta{\scriptstyle[0,0]} = 0$.
\item For any $(i,j)\in \mathbb{Z}^2$, $\theta{\scriptstyle[i,j]} = \theta{\scriptstyle[-i,-j]}$. 
\item For any $(\omega_1,\omega_2)\in[0,2\pi)^2$, $1-\sum_{( i,j)\in\mathbb{Z}^2}\theta{\scriptstyle[i,j]}\cos\left(i\omega_1+j\omega_2\right) > 0$.
\end{enumerate}
\end{defi}
Similarly, we define the set $\Theta^{+,\infty,\text{iso}}$ for the isotropic GMRFs on the lattices. As done in Section \ref{section_presentation} for  toroidal lattices, we now introduce the parametric 
 parameter sets. For any model $m\in\mathcal{M}_1$, the set 
$\Theta^{+,\infty}_m$ refers to the subset of matrices $\theta$ in $\Theta^{+,\infty}$ whose support is included in $m$.
Analogously, we define the parameter set $\Theta^{+,\infty,\text{iso}}_m$ corresponding to isotropic GMRFs.\\

We cannot directly extend the CLS empirical contrast $\gamma_{n,p_1,p_2}(.)$ defined in (\ref{definition_gamman}) in this new setting because we have to take the edge effect into account. Indeed, if we want to compute the conditional regression of ${\bf X}_i{\scriptstyle[j_1,j_2]}$, we have to observe \emph{all} its neighbors with respect to $m$, i.e. $\left\{{\bf X}_i{\scriptstyle[j_1+l_1,j_2+l_2]},\ (l_1,l_2)\in m\right\}$.  In this regard, we define the sublattice $\Lambda_m$  for any model $m\in\mathcal{M}_1$.
\begin{eqnarray*}
\Lambda_m:= \left\{(i_1,i_2)\in\Lambda\, ,\,\,(m+(i_1,i_2))\subset \Lambda\right\}\ ,  
\end{eqnarray*}
where $(m+(i,j))$ denotes the set $m$ of nodes  translated by $(i,j)$. For instance, if we consider the model $m_1$ with  four nearest neighbors, the edge effect size is  one and  $\Lambda_m$ contains all the nodes that do not lie on the border.
The model $m_3$ with 12 nearest neighbors yields an edge effect of size 2 and $\Lambda_m$ contains all the nodes in $\Lambda$, except those which are at a (euclidean) distance strictly smaller than $2$ from the border.

 For any model $m\in\mathcal{M}_1$, any $\theta'\in\Theta_m^{+,\infty}$, and any sublattice $\Lambda'\subset \Lambda_m$, we define  $\gamma_{n,p_1,p_2}^{\Lambda'}(.)$ as an analogous of $\gamma_{n,p_1,p_2}(.)$ except that it only relies on the conditional regression of the nodes in $\Lambda'$. 
\begin{eqnarray*}
\gamma^{\Lambda'}_{n,p_1,p_2}(\theta') & := & \frac{1}{n\text{Card}(\Lambda')}\sum_{i=1}^n\sum_{(j_1,j_2)\in \Lambda'} \bigg({\bf X}_i{\scriptstyle[j_1,j_2]} -
\sum_{(l_1,l_2)\in m} \theta'{\scriptstyle[l_1,l_2]}{\bf X}_i{\scriptstyle[j_1+l_1,j_2+l_2]}\bigg)^2 \ .
\end{eqnarray*}
Then, the \emph{CLS} estimators $\widehat{\theta}_m^{\Lambda'}$ and $\widehat{\theta}_m^{\Lambda',\text{iso}}$ are defined by 
\begin{eqnarray*}
 \widehat{\theta}_m^{\Lambda'} \in \arg \min_{\theta'\in \Theta^{+,\infty}_m} \gamma^{\Lambda'}_{n,p_1,p_2}\left(\theta'\right) \hspace{1.5cm}\text{and}\hspace{1.5cm} \widehat{\theta}_m^{\Lambda',\text{iso}} \in \arg \min_{\theta'\in \Theta^{+,\infty,\text{iso}}_m} \gamma^{\Lambda'}_{n,p_1,p_2}\left(\theta'\right)\ .
\end{eqnarray*}
Contrary to $\widehat{\theta}_m$, the estimator $\widehat{\theta}_m^{\Lambda_m}$ is not necessarily unique especially if the size of $\Lambda_m$ is smaller than $d_m$. Let us mention that it  is quite classical in the literature to remove nodes to  take edge effects or missing data into account (see e.g. \cite{guyon95} Sect.4.3). We cannot use anymore fast Fourier transform for computing the parametric estimator. Nevertheless, the estimators $\widehat{\theta}_m^{\Lambda'}$ are still computationally amenable, since they minimizes a quadratic function on the closed convex set $\Theta_m^{+,\infty}$.

Suppose we are given a subcollection $\mathcal{M}$ of $\mathcal{M}_1$. We note $\Lambda_{\mathcal{M}}$ the smallest sublattice  among the collection of lattices $\Lambda_m$ with $m\in\mathcal{M}$. In order to select the neighborhood $\widehat{m}$, we compute the estimators 
$\widehat{\theta}_m^{\Lambda_{\mathcal{M}}}$ and minimize the criteria $\gamma_{n,p_1,p_2}^{\Lambda_{\mathcal{M}}}(\widehat{\theta}_m^{\Lambda_{\mathcal{M}}})$ penalized by a quantity of the  order $d_m/(n\text{Card}(\Lambda_{\mathcal{M}}))$. We compute the quantities $\gamma_{n,p_1,p_2}^{\Lambda_{\mathcal{M}}}(\widehat{\theta}_m^{\Lambda_{\mathcal{M}}})$ instead of $\gamma_{n,p_1,p_2}^{\Lambda_m}(\widehat{\theta}_m^{\Lambda_m})$ since we want to compare the adequation of the models using the \emph{same} data set.

We now describe a data-driven model selection procedure for choosing the neighborhood. It is based on the slope heuristic developed in the previous section.  

\begin{algo}(Data-driven penalization  for non-toroidal lattice).\label{pente_nontorique}
 \begin{enumerate}
  \item Compute the selected model $\widehat{m}(N)$ as a function of $N>0$
$$\widehat{m}(N) \in \arg\min_{m\in\mathcal{M}}\left\{\gamma^{\Lambda_{\mathcal{M}}}_{n,p_1,p_2}(\widehat{\theta}^{\Lambda_{\mathcal{M}}}_m)+ N \frac{d_m}{n\text{Card}(\Lambda_{\mathcal{M}})}\right\}\ .$$ 
\item Find $\widehat{N}_{\text{\emph{min}}}>0$ such that the jump $d_{\widehat{m}\left(\left[\widehat{N}_{\text{\emph{min}}}\right]_-\right)}-d_{\widehat{m}\left(\left[\widehat{N}_{\text{\emph{min}}}\right]_+\right)}$ is maximal. 
\item Select the model  $\widehat{m}= \widehat{m}(2\widehat{N}_{\text{\emph{min}}})$.
\item Compute the estimator $\widehat{\theta}_{\widehat{m}}^{\Lambda_{\widehat{m}}}$.
 \end{enumerate}
\end{algo}
This procedure straightforwardly extends to the case of isotropic GMRFs estimation by replacing $\widehat{m}(N)$ by $\widehat{m}^{\text{iso}}(N)$ and $d_m$ by $d_{m}^{\text{iso}}$. For short, we write $\widetilde{\theta}$ (resp. $\widetilde{\theta}^{\text{iso}}$) for $\widehat{\theta}_{\widehat{m}}^{\Lambda_{\widehat{m}}}$ (resp. $\widehat{\theta}_{\widehat{m}}^{\Lambda_{\widehat{m}},\text{iso}}$).
As for Algorithm \ref{pente}, it is advised  to introduce huge models in the collection $\mathcal{M}$ in order to better detect the dimension jump. However, when the dimension of the models increases the size of $\Lambda_m$ decreases and the estimator $\widehat{\theta}^{\Lambda_m}_m$ may become unreliable. The method therefore requires a reasonable number of data. In practice, $\Lambda$ should not contain less than 100 nodes.

\section{Simulation study}\label{section_simulations}

In the first simulation experiment, we compare the efficiency of our procedure with penalized maximum likelihood methods when the field is a torus. In the second and third studies, we consider the estimation of a Gaussian field observed on a rectangle. The calculations are made with $R$~\cite{R}. Throughout these simulations, we only consider isotropic estimators.
	
\subsection{Isotropic GMRF on a torus}\label{simu-toriques}

First, we consider $X$ an isotropic GMRF on the torus $\Lambda$ of size $p=p_1=p_2=20$. There are therefore 400 points in the lattice. The number  of observations $n$ equals one and the conditional variance $\sigma^2$ is one. We introduce a radius $r:=\sqrt{17}$. Then, for any number $\phi>0$, we define the $p\times p$ matrix $\theta^\phi$ as:
\begin{eqnarray*}
\left\{\begin{array}{cccl}
 \theta^\phi{\scriptstyle[0,0]} & := & 0 \ ,&\\
\theta^\phi{\scriptstyle[i,j]} &:= & \phi &\text{if $|(i,j)|_t\leq r$ and $(i,j)\neq (0,0)$}\ ,\\
\theta^\phi{\scriptstyle[i,j]} &:= & 0 & \text{if $|(i,j)|_t> r$} \ .
\end{array}\right.
\end{eqnarray*}	
In practice, we set $\phi$ to $0$, $0.0125$, $0.015$, and $0.0175$. Observe that these choices constrain $\|\theta^\phi\|_1<1$. The matrix $\theta^\phi$ therefore belongs to the  set $\Theta_{m_{10}}^{+,\text{iso}}$ of dimension $10$ introduced in Definition \ref{definition_modele}.\\

{\bf First simulation experiment}. In Section \ref{section_resultats_theorique}, we have advocated the use of the estimator $\widetilde{\theta}$ instead of $\widetilde{\theta}_\rho$, although theoretical results are only available for $\widetilde{\theta}_\rho$ with $\rho<\infty$. We recall that  $\widetilde{\theta}=\widetilde{\theta}_{\rho}$ with $\rho=\infty$. We check in this simulation study that the performances of $\widetilde{\theta}$ and $\widetilde{\theta}_{\rho}$ with different values of $\rho$ are similar.

We consider the collection of neighborhoods $\mathcal{M}:=\left\{m_0,m_1,\ldots, m_{20}\right\}$ whose maximal dimension $d_{m_{20}}^{\text{iso}}$ is 21. The estimator $\widetilde{\theta}^{\text{iso}}$ is built using the CLS model selection procedure introduced in Algorithm \ref{pente}. The estimators $\widetilde{\theta}^{\text{iso}}_{\rho}$ are computed similarly, except that they are based on the parametric estimators $\widehat{\theta}^{\text{iso}}_{m,\rho}$ (Sect. \ref{section_resultats_theorique}) instead of $\widehat{\theta}^{\text{iso}}_{m}$.

The Gaussian field $X$  with $\phi=0.015$ is simulated by using the fast Fourier transform. 
The quality of the estimations is assessed by the prediction loss function $l(.,.)$ defined in (\ref{definition_perte}).  The experiments are repeated $1000$ times. For $\rho=2,\ 4,\ 8$, we evaluate the risks $\mathbb{E}_{\theta^{\phi}}[l(\widetilde{\theta}^{\text{iso}},\theta^{\phi})]$ and $\mathbb{E}_{\theta^{\phi}}[l(\widetilde{\theta}^{\text{iso}}_{\rho},\theta^{\phi})]$
as well as the corresponding empirical $95\%$ confidence intervals  by a Monte-Carlo method. We also estimate the risks of $\widehat{\theta}^{\text{iso}}_m$ and  $\widehat{\theta}^{\text{iso}}_{m,\rho}$ for each model $m\in\mathcal{M}$. It then allows to evaluate the  oracle risks $\mathbb{E}_{\theta^{\phi}}[l(\widehat{\theta}^{\text{iso}}_{m^*,\rho},\theta^{\phi})]$
and the risk ratios $\mathbb{E}_{\theta^{\phi}}[l(\widetilde{\theta}^{\text{iso}}_{\rho},\theta^{\phi})]/\mathbb{E}_{\theta^{\phi}}[l(\widehat{\theta}^{\text{iso}}_{m^*,\rho},\theta^{\phi})]$. The risk ratio measures how well the selected model $\widehat{m}^{\text{iso}}$ performs in comparison to the ``best'' model $m^*$. Moreover, the risk ratio  roughly illustrates the oracle type inequality presented in Theorem \ref{mainthrm}. Indeed, the infimum $\inf_{m\in\mathcal{M}}[l(\theta_{m,\rho},\theta)+\pen(m)]$ in (14) is a good measure of the risk $\mathbb{E}_{\theta^{\phi}}[l(\widehat{\theta}^{\text{iso}}_{m^*,\rho},\theta^{\phi})]$ as explained in \cite{verzelen_gmrf_theorie} Sect.4. The results are given in Table \ref{resultat0}. They corroborate that the estimators $\widetilde{\theta}^{\text{iso}}$ and $\widetilde{\theta}^{\text{iso}}_\rho$ perform similarly. Moreover, the risk ratios $\mathbb{E}_{\theta^{\phi}}[l(\widetilde{\theta}^{\text{iso}}_{\rho},\theta^{\phi})]/\mathbb{E}_{\theta^{\phi}}[l(\widehat{\theta}^{\text{iso}}_{m^*,\rho},\theta^{\phi})]$ correspond to the ratios

\begin{Table}[h]	
\caption{First simulation study. Estimates and $95\%$ confidence intervals of the risks $\mathbb{E}_{\theta^{\phi}}[l(\widetilde{\theta}^{\text{iso}},\theta^{\phi})]$, $\mathbb{E}_{\theta^{\phi}}[l(\widetilde{\theta}^{\text{iso}}_{\rho},\theta^{\phi})]$, and of the ratios $\mathbb{E}_{\theta^{\phi}}[l(\widetilde{\theta}^{\text{iso}},\theta^{\phi})]/\mathbb{E}_{\theta^{\phi}}[l(\widehat{\theta}^{\text{iso}}_{m^*},\theta^{\phi})]$ and $\mathbb{E}_{\theta^{\phi}}[l(\widetilde{\theta}^{\text{iso}}_{\rho},\theta^{\phi})]/\mathbb{E}_{\theta^{\phi}}[l(\widehat{\theta}^{\text{iso}}_{m^*,\rho},\theta^{\phi})]$ with $\phi=0.015$ and $\rho=2,\ 4,\ 8$.   \label{resultat0}}
\begin{center}
\begin{tabular}{c|c|c|c|c}
$\rho$  & 2   & 4 & 8 & $\infty$\\
\hline
$\mathbb{E}_{\theta^{\phi}}[l(\widetilde{\theta}^{\text{iso}}_{\rho},\theta^{\phi})]\times 10^2$& $4.1\pm0.1$&   $4.2\pm 0.2$ &$4.2\pm 0.1$ & $4.2 \pm 0.3$ \\
$\mathbb{E}_{\theta^{\phi}}[l(\widetilde{\theta}^{\text{iso}}_{\rho},\theta^{\phi})]/\mathbb{E}_{\theta^{\phi}}[l(\widehat{\theta}^{\text{iso}}_{m^*,\rho},\theta^{\phi})]$ &$1.3\pm 0.1$  & $1.3\pm 0.1$ & $1.3\pm 0.1$ &  $1.3\pm 0.2$
\end{tabular}
\end{center}
\end{Table}

\vspace{0.5cm}

{\bf Second simulation experiment}.
We compare the efficiency of the method with two alternative model selection procedures. For each of them, we use the collection $\mathcal{M}$ as in the previous experiment. The two alternative procedures are based on likelihood maximization.
 In this regard, we first define the parametric maximum likelihood estimator $\widehat{\theta}^{\text{mle}}_m$ for any model $m\in\mathcal{M}$,
\begin{eqnarray*}
\left(\widehat{\theta}^{\text{mle}}_m,\widehat{\sigma}^{\text{mle}}_m\right):= \arg\min_{\theta'\in\Theta^{+,\text{iso}}_m,\sigma'} -\mathcal{L}_p(\theta',\sigma',{\bf X})\ ,
\end{eqnarray*}
where $\mathcal{L}_p(\theta',{\bf X})$ stands for  the $\log$-likelihood at the parameter $\theta'$. 
We then select a model $m$ applying either an  $\AIC$-type criterion \cite{akaike73} or a $\BIC$-type criterion \cite{BIC}:
\begin{eqnarray*}
 \widehat{m}^{\text{AIC}} &:=& \arg \min_{m\in\mathcal{M}}\left\{-2\mathcal{L}_p(\widehat{\theta}^{\text{mle}}_m,\widehat{\sigma}^{\text{mle}}_m,{\bf X})+ 2d_m^{\text{iso}}\right\}\ , \\
\widehat{m}^{\text{BIC}} &:=& \arg \min_{m\in\mathcal{M}}\left\{-2\mathcal{L}_p(\widehat{\theta}^{\text{mle}}_m,\widehat{\sigma}^{\text{mle}}_m,{\bf X})+ \log(p^2)d_m^{\text{iso}}\right\}\ .
\end{eqnarray*}
For short, we write $\widehat{\theta}^{\text{AIC}}$ and $\widehat{\theta}^{\text{BIC}}$ for the two obtained estimators $\widehat{\theta}^{\text{mle}}_{\widehat{m}^{\text{AIC}}}$ and $\widehat{\theta}^{\text{mle}}_{\widehat{m}^{\text{BIC}}}$. 
Although $\AIC$ and $\BIC$ procedures are not justified in this setting, we still apply them as they are widely used in many frameworks. Their computation is performed efficiently using the fast Fourier transform described  in Section \ref{section_computationelle}. \\

The experiments are repeated $1000$ times. The Gaussian field is simulated using the fast Fourier transform. 
The quality of the estimations is assessed by the prediction loss function $l(.,.)$.
For any $\phi$ and any of these three estimators, 
 we evaluate the risks $\mathbb{E}_{\theta^{\phi}}[l(\widehat{\theta}^{\text{AIC}},\theta^{\phi})] $, $\mathbb{E}_{\theta^{\phi}}[l(\widehat{\theta}^{\text{BIC}},\theta^{\phi})]$, and $\mathbb{E}_{\theta^{\phi}}[l(\widetilde{\theta}^{\text{iso}},\theta^{\phi})]$ as well as the corresponding empirical $95\%$ confidence intervals  by a Monte-Carlo method. We also estimate  the risk ratios $\mathbb{E}_{\theta^{\phi}}[l(\widetilde{\theta}^{\text{iso}},\theta^{\phi})]/\mathbb{E}_{\theta^{\phi}}[l(\widehat{\theta}^{\text{iso}}_{m^*},\theta^{\phi})]$
 The results  are given in Table \ref{resultat1}.\vspace{0.5cm}

\begin{Table}[h]	
\caption{Second simulation study. Estimates and $95\%$ confidence intervals of the risks $\mathbb{E}_{\theta^{\phi}}[l(\widehat{\theta}^{\text{AIC}},\theta^{\phi})]$,  $\mathbb{E}_{\theta^{\phi}}[l(\widehat{\theta}^{\text{BIC}},\theta^{\phi})]$, and $\mathbb{E}_{\theta^{\phi}}[l(\widetilde{\theta}^{\text{iso}},\theta^{\phi})]$ and of the ratio $\mathbb{E}_{\theta^{\phi}}[l(\widetilde{\theta}^{\text{iso}},\theta^{\phi})]/\mathbb{E}_{\theta^{\phi}}[l(\widehat{\theta}^{\text{iso}}_{m^*},\theta^{\phi})]$.   \label{resultat1}}
\begin{center}
\begin{tabular}{c|c|c|c|c}
$\phi\times 10^2$  & 0   & 1.25 & 1.5 & 1.75\\
\hline
$\mathbb{E}_{\theta^{\phi}}[l(\widehat{\theta}^{\text{AIC}},\theta^{\phi})]\times 10^2$ &$1.2\pm 0.2$ &$3.1\pm 0.2$ &$4.3\pm 0.2$ & $6.4\pm 0.2$\\
$\mathbb{E}_{\theta^{\phi}}[l(\widehat{\theta}^{\text{BIC}},\theta^{\phi})]\times 10^2$& $0.01\pm 0.01$  & $1.9\pm 0.1$ & $3.7\pm 0.1$&$9.7 \pm0.3 $\\
$\mathbb{E}_{\theta^{\phi}}[l(\widetilde{\theta}^{\text{iso}},\theta^{\phi})]\times 10^2$& $1.6\pm0.2$&   $3.2\pm 0.2$ &$4.2\pm 0.1$ & $7.2 \pm 0.3$ \\
$\mathbb{E}_{\theta^{\phi}}[l(\widetilde{\theta}^{\text{iso}},\theta^{\phi})]/\mathbb{E}_{\theta^{\phi}}[l(\widehat{\theta}^{\text{iso}}_{m^*},\theta^{\phi})]$ & $+\infty $ & $1.9\pm 0.7$ & $1.3\pm 0.2$ & $1.5\pm 0.3$
\end{tabular}

\end{center}
\end{Table}

The $\BIC$ criterion outperforms the other procedures when $\phi=0$, $0.0125$, or $0.015$ but behaves bad for a large $\phi$. Indeed, the $\BIC$ criterion has a tendency to overpenalize the models. For the two first values  of $\phi$ the oracle model in $\mathcal{M}$ is $m_0$. Hence, overpenalizing  increases the performance of estimation in this case. However, when $\phi$ increases, the dimension of the oracle model is larger and $\BIC$ therefore selects too small models.

In contrast, $\AIC$ and the CLS estimator exhibit similar behaviors.   If we forget the case $\phi=0$ for which the oracle risk is 0, the risk of $\widetilde{\theta}^{\text{iso}}$ is close to the risk of the oracle model (the ratio is close to one). Hence, the neighborhood choice for $\widetilde{\theta}^{\text{iso}}$ is almost optimal. 

In conclusion, $\widetilde{\theta}^{\text{iso}}$ or $\widehat{\theta}^{\text{AIC}}$ both exhibit good performances
for estimating the distribution of a regular Gaussian field on a torus. The strength of our neighborhood selection procedure lies in the fact it easily generalizes to non-toroidal lattices as illustrated in the next section.

\subsection{Isotropic Gaussian fields on $\mathbb{Z}^2$}\label{simu-non-toriques}

{\bf First simulation experiment}.
We now consider $X$ an isotropic Gaussian field defined on $\mathbb{Z}^2$ but only observed on a square $\Lambda$ of sizes $p=p_1=p_2=20$ or  $p=p_1=p_2=100$. This corresponds to the setting described in Section \ref{section_extension_non_torique}. The variance of $X{\scriptstyle[0,0]}$ is set to one and the distribution of the field is therefore uniquely defined by its correlation function $\rho(k,l):=\corr(X{\scriptstyle[k,l]},X{\scriptstyle[0,0]})$. 
Again, the number of replications $n$ is chosen to be one.
In the first experiment, we use four classical correlation functions: exponential, spherical, circular, and Mat\'ern (e.g. \cite{cressie} Sect.2.3.1 and \cite{matern86}).
\begin{eqnarray*}
\text{Exponential:} & \rho(k,l) =&  \exp\left(-\frac{d(k,l)}{r}\right) \\
\text{Circular:} & \rho(k,l) = &\left\{\begin{array}{ccl}
1-\frac{2}{\pi}\left[\frac{d(k,l)}{r}\sqrt{1-\left(\frac{d(k,l)}{r}\right)^2}+\sin^{-1} \left(\sqrt{\frac{d(k,l)}{r}}\right)\right] & \text{if}& d(k,l)\leq r\\
0 & \text{else} & \end{array}\right.\\
\text{Spherical:} & \rho(k,l) = &\left\{\begin{array}{ccl}
1-1.5\frac{d(k,l)}{r}+0.5 \left(\frac{d(k,l)}{r}\right)^3 & \text{if}& d(k,l)\leq r\\
0 & \text{else} & 
\end{array}
\right.\\
\text{Mat\'ern:} & \rho(k,l) =&	 \frac{1}{2^{\kappa-1}\Gamma(\kappa)}\left(\frac{d(k,l)}{r}\right)^{\kappa}\mathcal{K}_{\kappa}\left(\frac{d(k,l)}{r}\right)\ , 
\end{eqnarray*}
where $d(k,l)$ denotes the euclidean distance from $(k,l)$ to $(0,0)$ and $\mathcal{K}_{\kappa}(.)$ is the modified Bessel function of order $\kappa$. In a nutshell, the parameter $r$ represents the range of correlation, whereas $\kappa$ may be regarded as a smoothness parameter for the Mat\'ern function. In this simulation experiment, we set $r$ to $3$. When considering the Mat\'ern model, we take $\kappa$ equal to $0.05$, $0.25$, $0.5$, $1$, $2$, and $4$.

The Gaussian fields are simulated using the function \emph{GaussRF} in the library \emph{RandomFields}~\cite{RF}.  For each of experiments, 
we compute the estimator $\widetilde{\theta}^{\text{iso}}$ based on Algorithm \ref{pente_nontorique} with the collection $\mathcal{M}:=\{m\in\mathcal{M}_1\, , \,d_m^{\text{iso}}\leq 18\}$. Since the lattice $\Lambda$ is not a torus, methods based on likelihood maximization exhibit a prohibitive computational burden. Consequently, we do not use MLE in this experiment. 
We shall compare the efficiency of  $\widetilde{\theta}^{\text{iso}}$ with a variogram-based  estimation method. 

We recall that the linear combination $\sum_{(i,j)\in\Lambda\setminus\{(0,0)\}}\theta{\scriptstyle[i,j]}X{\scriptstyle[i,j]}$ is the kriging predictor of $X{\scriptstyle[0,0]}$ given the remaining variables (Equation (\ref{regression_conditionnelle})).  A natural method to estimate $\theta$ in this spatial setting amounts to estimating the variogram of the observed Gaussian field and then  performing ordinary kriging at the node $(0,0)$. More precisely,  we 
first estimate the empirical variogram by applying the modulus estimator of Hawkes and Cressie (e.g. \cite{cressie} Eq.(2.2.8)) to the observed field of $400$ points. Afterwards, we fit this empirical variogram to a variogram model using the reweighted least-squares suggested by Cressie \cite{cressie85}. This procedure therefore requires the choice of a particular variogram model. In the first simulation study, we choose \emph{the model} that has generated the data. Observe that this method is \emph{not} adaptive since it requires the knowledge of the variogram model.
In practice, we use Library \emph{geoR}~\cite{geoR} implemented in \emph{R}~\cite{R}  to estimate the parameters $r$, $\var(X{\scriptstyle[0,0]})$ and eventually $\kappa$ of the variogram model. Then, we compute the estimator $\widehat{\theta}^V$ by performing ordinary kriging at the center node of $\Lambda$. For each of these estimations, we assume that the variogram model is known.
For computational reasons, we use a kriging neighborhood of size $11\times 11$ that contains 120 points. Previous simulations have indicated that this neighborhood choice does not decrease the precision of the estimation. For the Mat\`ern model with $\kappa=2$ and $4$, the covariance is almost singular. There are sometimes inversion difficulties and we therefore use kriging neighborhood of respective size $7\times 7$ and $3\times 3$.

We again assess the performances of the procedures using the loss $l(.,.)$. Even if this loss is defined in (\ref{definition_perte}) for a torus, the alternative definition (\ref{definition_perte_alternative}) clearly extends to this non-toroidal setting. Consequently, the loss $l(\widehat{\theta},\theta)$ measures the difference between the prediction error of $X{\scriptstyle [0,0]}$ when using $\sum_{(i,j)\in \Lambda \setminus\{(0,0)\}}\widehat{\theta}{\scriptstyle [i,j]} X{\scriptstyle [i,j]}$ and the prediction error of $X{\scriptstyle [0,0]}$ when using the best predictor $\mathbb{E}[X{\scriptstyle [0,0]}|(X{\scriptstyle [i,j]})_{(i,j)\in\Lambda\setminus\{(0,0)\}}]$. In other words, $l(\widehat{\theta},\theta)$ is the difference of the kriging error made with the estimated parameters $\widehat{\theta}$ and the kriging error made with the true parameter $\theta$. \\

The experiments are repeated $1000$ times.
For any of the four correlation models previously mentioned, we evaluate the risks $\mathbb{E}_{\theta}[l(\widetilde{\theta}^{\text{iso}},\theta)]$ and $\mathbb{E}_{\theta}[l(\widehat{\theta}^V,\theta)]$ by Monte-Carlo. In order to assess the efficiency of the selection procedure, we also evaluate the risk ratio 
\begin{eqnarray*}
Risk.ratio= \frac{\mathbb{E}_{\theta}[l(\widehat{\theta}^{\Lambda_{\mathcal{M}},\text{iso}}_{\widehat{m}},\theta)]}{\mathbb{E}_{\theta}[l(\widehat{\theta}^{\Lambda_{\mathcal{M}},\text{iso}}_{m^*},\theta)]} \ .
\end{eqnarray*}
 As in Section \ref{simu-toriques}, the oracle risk $\mathbb{E}[l(\widehat{\theta}^{\Lambda_{\mathcal{M}},\text{iso}}_{m^*},\theta)]$ is evaluated by taking the minimum of the evaluations of the risks $\mathbb{E}[l(\widehat{\theta}^{\Lambda_{\mathcal{M}},\text{iso}}_{m},\theta)]$ over all models $m\in\mathcal{M}$. Results of the simulation experiment are given in Table \ref{resultat2} and \ref{resultat3}.

Observe that none of the fields considered in this study are GMRFs. Here, the GMRF models should only be viewed as a collection of approximation sets of the true distribution. This simulation experiment is in the spirit of  Rue and Tjelmeland's study \cite{rue02}. However, there are some major differences. Contrary to them, we perform estimation and not only approximation. Moreover, our lattice is not a torus. Finally, we use our prediction loss $l(.,.)$ to assess the performance, whereas they compare the correlation functions.

\begin{Table}[h]
\caption{Estimates and $95\%$ confidence intervals of the risks $\mathbb{E}_{\theta}[l(\widehat{\theta}^V,\theta)]$ and  $\mathbb{E}_{\theta}[l(\widetilde{\theta}^{\text{iso}},\theta)]$ and of $Risk.ratio$ for the exponential, circular and spherical models with $p=20$.   \label{resultat2}}
\begin{center}

\begin{tabular}{c|c|c|c|c}
 Model & Exponential & Circular & Spherical  \\
\hline
$\mathbb{E}_{\theta}[l(\widehat{\theta}^V,\theta]\times 10^2$& $0.08\pm 0.01$ & $9.1\pm 0.5$ & $2.9\pm 0.1 $  \\
$\mathbb{E}_{\theta}[l(\widetilde{\theta}^{\text{iso}},\theta)]\times 10^2$& $1.08\pm 0.01$ & $6.5\pm 0.1$ & $3.4\pm 0.1$ \\
$Risk.ratio$ & $3.6\pm 0.4$ &$1.4\pm0.1$ &$1.6\pm 0.1$ 

\end{tabular}

\end{center}
\end{Table}

\begin{Table}[h]
\caption{Estimates and $95\%$ confidence intervals of the risks $\mathbb{E}_{\theta}[l(\widehat{\theta}^V,\theta)]$ and  $\mathbb{E}_{\theta}[l(\widetilde{\theta}^{\text{iso}},\theta)]$ and of $Risk.ratio$  for Mat\'ern model with $p=100$.   \label{resultat3}}
\begin{center}
\begin{tabular}{c|c|c|c|c|c|}	
 $\kappa$  & 0.05 & 0.25 & 0.5 & 1\\
\hline
$\mathbb{E}_{\theta}[l(\widehat{\theta}^V,\theta)]\times 10^3$ &$91.8\pm 0.7$ & $80.0\pm 0.2$ & $18.0\pm 0.1$ & $2.5\pm 0.1$	 \\
$\mathbb{E}_{\theta}[l(\widetilde{\theta}^{\text{iso}},\theta)]\times 10^3$ &$2.24\pm 0.01$ & $0.62\pm 0.01$ &$0.33\pm 0.01$& $0.08\pm 0.01$\\
$Risk.ratio$ & $1.3\pm 0.1$ & $1.7	\pm 0.2$ & $1.5\pm 0.2$ & $1.3\pm 0.1$
\end{tabular}

\vspace{0.4cm}

\begin{tabular}{c|c|c|}	
$\kappa$  &2 & 4\\
\hline $\mathbb{E}_{\theta}[l(\widehat{\theta}^V,\theta)]\times 10^4$ & $6.3\pm 1.1$ & $0.011\pm 0.001$\\
$\mathbb{E}_{\theta}[l(\widetilde{\theta}^{\text{iso}},\theta)]\times 10^4$& $1.9\pm 0.1$ & $0.17\pm 0.01$\\
$Risk.ratio$ & $2.6\pm 0.2$ & $1.1\pm 0.1$
\end{tabular}
\end{center}
\end{Table}

\emph{Comments on Tables \ref{resultat2} and \ref{resultat3}}. In both tables, the ratio $\mathbb{E}_{\theta}[l(\widehat{\theta}^{\Lambda_{\mathcal{M}},\text{iso}}_{\widehat{m}},\theta)]/\mathbb{E}_{\theta}[l(\widehat{\theta}^{\Lambda_{\mathcal{M}},\text{iso}}_{m^*},\theta)]$ stays close to one. Hence, the model selection is almost optimal from an efficiency point of view. In most of the cases, the estimator $\widetilde{\theta}^{\text{iso}}$ outperforms the estimator $\widehat{\theta}^V$ based on geostatistical methods. This is particularly striking for the Mat\'ern correlation model because in that case the computation of  $\widehat{\theta}^V$ requires the estimation of the additional parameter $\kappa$. Indeed, let us recall that the exponential model and the Mat\'ern model with $\kappa=0.5$ are equivalent. For $\kappa=0.5$, the risk of $\widehat{\theta}^V$ is $100$ times higher when $\kappa$ has to be estimated than when $\kappa$ is known.\\

{\bf Second simulation experiment}. The kriging estimator $\widehat{\theta}^V$ requires the knowledge or the choice of a correlation model. In the second simulation experiment, the correlation of $X$ is the Mat\`ern function with range $r=3$ and $\kappa=0.05$. The size $p$ of the lattice is chosen to be $100$. We now estimate $\theta$ using \emph{different} variogram models, namely the exponential, the circular, the spherical and the Mat\`ern model. The estimator $\widetilde{\theta}^{\text{iso}}$ for such a field was already considered in Table \ref{resultat3}. The experiment is repeated $1000$ times.

\begin{Table}[h]
\caption{Estimates and $95\%$ confidence intervals of the risks $\mathbb{E}_{\theta}[l(\widehat{\theta}^V,\theta)]$  for Mat\'ern model with $\kappa=0.05$ when using the exponential, circular, spherical, and Mat\`ern models with $p=100$.   \label{resultat_mauvais_variogramme}}
\begin{center}
\begin{tabular}{c|c|c|c|c|c}	
 Model  & Exponential  & Circular & Spherical & Mat\`ern\\
\hline
$\mathbb{E}_{\theta}[l(\widehat{\theta}^V,\theta)]\times 10^3$ & $48.3\pm 0.4$ & $461\pm 16$ & $293\pm 7$ & $91.8 \pm 0.7$ 
\end{tabular}
\end{center}

\end{Table}

\emph{Comments on Table \ref{resultat_mauvais_variogramme}}. One observes that circular and spherical models yield worse performances than Mat\`ern model. In contrast, the exponential model behaves better. The choice of the variogram model therefore seems critical to get good performances. The model selection estimator $\widetilde{\theta}^{\text{iso}}$ (Table \ref{resultat3}) exhibits a smaller risk than the exponential model.\\

\subsection{Anisotropic Gaussian fields on $\mathbb{Z}^2$}\label{simu-anisotropic}

We still consider $X$ a Gaussian field observed on a square  $\Lambda$ of size $100\times 100$. Contrary to the previous study, the field is not assumed to be isotropic. To model the geometric anisotropy, we suppose that $X$ is an isotropic  field on a deformed lattice $\Lambda'$. The transformation consists in multiplying the original coordinates by a rotation $R$ and a shrinking matrix $T$. For the sake of simplicity, we take the identity for $R$. The shrinking matrix $T$ is defined by the anisotropy ratio (\emph{Ani.ratio}). It  corresponds to the ratio between the directions with smaller and greater continuity in the field $X$, i.e the ratio between maximum and minimum ranges. In this experiment,   $X$ follows a  Mat\`ern correlation  with range $r=3$, $\kappa=0.05$, 0.25, 0.5, 1, 2, and 4 and \emph{Ani.ratio}=2 or 5. We compute the anisotropic estimator $\widetilde{\theta}$ based on Algorithm \ref{pente_nontorique} with the collection $\mathcal{M}:=\left\{m\in\mathcal{M}_1,d_m\leq 28\right\}$. As a benchmark, we also compute the variogram-based estimator $\widehat{\theta}^V$ based on the Mat\`ern model. In order to compute $\widehat{\theta}^V$, we assume that we \emph{know} the anisotropy ratio and the anisotropy directions. Observe that the estimator $\widetilde{\theta}$ does not require any assumption on the form of anisotropy, while $\widehat{\theta}^V$ uses the geometric parameters of the anisotropy.\\

The experiments are repeated $1000$ times. We evaluate the risks $\mathbb{E}_{\theta}[l(\widehat{\theta}^V,\theta)]$ and  $\mathbb{E}_{\theta}[l(\widetilde{\theta},\theta)]$ and the risk ratio defined by 
\begin{eqnarray*}
 Risk.ratio = \frac{\mathbb{E}_{\theta}[l(\widehat{\theta}^{\Lambda_{\mathcal{M}}}_{\widehat{m}},\theta)]}{\mathbb{E}_{\theta}[l(\widehat{\theta}^{\Lambda_{\mathcal{M}}}_{m^*},\theta)]}\ .
\end{eqnarray*}

\begin{Table}[h]
\caption{Estimates and $95\%$ confidence intervals of the risks $\mathbb{E}_{\theta}[l(\widehat{\theta}^V,\theta)]$ and  $\mathbb{E}_{\theta}[l(\widetilde{\theta},\theta)]$ and of $Risk.ratio$ for Mat\'ern model and \emph{Ani.ratio}$=2$.   \label{resultat-ani1}}
\begin{center}
\begin{tabular}{c|c|c|c|c|c}	
 $\kappa$  & 0.05 & 0.25 & 0.5 & 1\\
\hline
$\mathbb{E}_{\theta}[l(\widehat{\theta}^V,\theta)]\times 10^2$ &$15.8\pm 0.1$ & $13.9\pm 0.1$ & $3.3\pm 0.1$ & $0.30\pm 0.01$\\	
$\mathbb{E}_{\theta}[l(\widetilde{\theta},\theta)]\times 10^2$ &$0.65\pm 0.01$ & $0.20\pm 0.01$ &$0.089\pm 0.001$& $0.17\pm 0.01$\\
$Risk.ratio$ & $1.2\pm 0.1$ & $1.1\pm 0.1$ & $1.1\pm 0.1$ & $1.7\pm 0.2$ 
\end{tabular}
\vspace{0.4cm}

\begin{tabular}{c|c|c|}	
$\kappa$  &2 & 4\\
\hline $\mathbb{E}_{\theta}[l(\widehat{\theta}^V,\theta)]\times 10^4	$ & $9.8\pm 0.1 $ & $ 0.020\pm 0.001$\\
$\mathbb{E}_{\theta}[l(\widetilde{\theta}^{\text{iso}},\theta)]\times 10^4$& $ 45.0\pm 0.1$ & $4.3\pm 0.1 $\\
$Risk.ratio$ & $ 2.9\pm 0.2$ & $22.3\pm 1.7 $
\end{tabular}

\end{center}
\end{Table}

\begin{Table}[h]
\caption{Estimates and $95\%$ confidence intervals of the risks $\mathbb{E}_{\theta}[l(\widehat{\theta}^V,\theta)]$ and $\mathbb{E}_{\theta}[l(\widetilde{\theta},\theta)]$ and of $Risk.ratio$ for Mat\'ern model and \emph{Ani.ratio}$=5$.   \label{resultat-ani2}}
\begin{center}
\begin{tabular}{c|c|c|c|c|c}	
 $\kappa$  & 0.05 & 0.25 & 0.5 & 1\\
\hline
$\mathbb{E}_{\theta}[l(\widehat{\theta}^V,\theta)]\times 10^2$ &$11.2\pm 0.1$ & $14.9\pm 0.1$ & $3.7\pm 0.1$ & $2.9\pm 0.1$\\	
$\mathbb{E}_{\theta}[l(\widetilde{\theta},\theta)]\times 10^2$ &$0.66\pm 0.1$ & $0.40\pm 0.01$ &$0.081\pm 0.001$& $0.14\pm 0.01$ \\
$Risk.ratio$ & $1.1\pm 0.1$ & $1.1\pm 0.1$ & $1.2\pm 0.1$ & $3.4\pm 0.8$ 
\end{tabular}

\vspace{0.4cm}

\begin{tabular}{c|c|c|}	
$\kappa$  &2 & 4\\
\hline $\mathbb{E}_{\theta}[l(\widehat{\theta}^V,\theta)]\times 10^4	$ & $30.6\pm 0.1 $ & $ 0.22\pm 0.01$\\
$\mathbb{E}_{\theta}[l(\widetilde{\theta}^{\text{iso}},\theta)]\times 10^4$& $ 38.0\pm 0.1$ & $39.6\pm 0.1 $\\
$Risk.ratio$ & $ 2.1\pm 0.1$ & $9.0\pm 1.4 $
\end{tabular}
\end{center}
\end{Table}

\emph{Comments on Tables \ref{resultat-ani1} and \ref{resultat-ani2}}. Except for the cases $\kappa=2,4$, the estimator $\widetilde{\theta}$ performs better than the variogram-based estimator $\widehat{\theta}^V$, although $\widehat{\theta}^V$ uses the true anisotropy parameters. For $\kappa=4$, the neighborhood selection is no performed efficiently (the risk ratio is large).

\section{Discussion}\label{section_conclusion}

In this paper, we have extended a neighborhood selection procedure introduced in \cite{verzelen_gmrf_theorie}. On the one hand, an algorithm is provided for tuning the penalty in practice. On the other hand, the new method also handles non-toroidal lattices. The computational complexity remains reasonable even when the size of the lattice is large.

In the case of stationary fields on a torus, our neighborhood selection procedure exhibits a computational burden and statistical performances analogous to the AIC procedure. Even if AIC has not been analyzed from an efficiency point of view, this suggests that AIC may achieve an oracle inequality in this setting. Moreover, we have empirically checked that $\widetilde{\theta}$ performs almost as well as the oracle model since the oracle ratio $\mathbb{E}[l(\widetilde{\theta},\theta)]/\mathbb{E}[l(\widehat{\theta}_{m^*},\theta)]$ remains close to one.

The strength of this neighborhood selection procedure lies in the fact it easily extends to non-toroidal lattices. We have illustrated that our method often outperforms variogram-based estimation methods in terms of the mean-squared prediction error. Moreover, the procedure behaves almost as well as the oracle. In contrast, variogram-based procedures may perform well for some covariances structure but also yield terrible results for other covariance structures. These results illustrate the \emph{adaptivity} of the neighborhood selection procedure.\\

In many statistical applications, Gaussian fields (or Gaussian Markov random fields) are not directly observed. For instance, Aykroyd \cite{aykroyd} or Dass and Nair \cite{dass03} use compound Gaussian Markov random fields to account for non stationarity and steep variations. The wavelet transform has emerged as a powerful tool in image analysis. The wavelet coefficients of an image are sometimes modeled using hidden Markov models~\cite{crouse98,portilla}. More generally, the success of the GMRFs is mainly due to the use of hierarchical models involving latent GMRFs~\cite{inla}. The study and the implementation of our penalization strategy for selecting the complexity of latent Markov models is an interesting direction of research.

\section{Proofs}\label{section_proofs}

Let us introduce some notations that shall be used throughout the proofs. For any $1\leq k\leq n$, the vector ${\bf X}^v_k$ denotes the vectorialized version of the $k$-th sample of $X$. Moreover, ${\bf X^v}$ is the matrix of size $p_1p_2\times n$ of the $n$ realisations of the vector ${\bf X}^v_k$. Throughout these proofs, $L,L_1, L_2$ denote constants that may vary from line to line. The notation $L(.)$ specifies the dependency on some quantities. Finally, the $\gamma(.)$ function stands for an infinite sampled version of the CLS criterion $\gamma_{n,p_1,p_2}(.)$:
$\gamma(.):=\mathbb{E}[\gamma_{n,p_1,p_2}(.)]$.

\subsection{Proof of Lemma \ref{lemme_transformee_fourier}}

Let us provide an alternative expression of $\gamma_{n,p_1,p_2}(\theta')$ in term of the factor $C(\theta')$ 
and the empirical covariance matrix $\overline{{\bf X^v X^{v*}}}$.
\begin{eqnarray}\label{definition_gamman_alternative}
\gamma_{n,p_1,p_2}(\theta')& = & \frac{1}{np_1p_2}tr\left[(I_{p_1p_2}-C(\theta'))\overline{{\bf X^v X^{v*}}}(I_{p_1p_2}-C(\theta'))\right] \ . 
\end{eqnarray}
This is justified in \cite{verzelen_gmrf_theorie} Sect.2.2.

\begin{lemma}\label{codiagonalisation}
There exists an orthogonal matrix $P$ which simultaneously diagonalizes every $p_1p_2\times
p_1p_2$ symmetric block circulant matrices with $p_2\times p_2$ blocks. Let $\theta$
be a  matrix of size $p_1\times p_2$ such that $C(\theta)$ is symmetric.
The matrix $D(\theta)= P^*C(\theta)P$ is diagonal and satisfies
\begin{eqnarray}\label{forme_valeurs_propres}
D(\theta) {\scriptstyle[(i-1)p_2+j,(i-1)p_2+j]} = \sum_{k=1}^{p_1}\sum_{l=1}^{p_2}\theta{\scriptstyle[k,l]}\cos\left[2\pi(ki/p_1+lj/p_2)\right]\ ,
\end{eqnarray}
for any $1\leq i\leq p_1$ and $1\leq j\leq p_2$.
\end{lemma}
This lemma is proved as in \cite{rue} Sect.2.6.2 to the price of a slight modification that takes into account the fact that $P$ is orthogonal and not unitary. The difference comes from the fact that contrary to Rue and Held we also assume that $C(\theta)$ is symmetric.
Lemma \ref{codiagonalisation} states that all symmetric block circulant matrices are simultaneously diagonalizable. Observe that for any $1\leq i\leq p_1$ and $1\leq j\leq p_2$, it holds that $D(\theta) {\scriptstyle[(i-1)p_2+j,(i-1)p_2+j]}=\lambda{\scriptstyle[i,j]}(\theta)$ since $\theta{\scriptstyle[k,l]}=\theta{\scriptstyle[p_1-k,p_2-l]}$. Hence, Expression (\ref{definition_gamman_alternative}) becomes
\begin{eqnarray*}
\gamma_{n,p_1,p_2}(\theta')& = & \frac{1}{np_1p_2}\bigg\{\sum_{i=1}^{p_1}\sum_{j=1}^{p_2}\left[1-\lambda{\scriptstyle[i,j]}(\theta)\right]^2\bigg[\sum_{k=1}^n \left[ P^*{\bf X}_{k}^{v}({\bf X}_{k}^{v})^*P\right]{\scriptstyle[(i-1)p_2+j,(i-1)p_2+j]} \bigg]   \bigg\}\ ,
\end{eqnarray*}
where ${\bf X}_{k}^{v}$ is the vectorialized version of the $k$-th observation of the field $X$. Straightforward computations allow us to prove that the quantities 
\begin{eqnarray*}	
\left(P^*{\bf X}_{k}^{v}({\bf X}_{k}^{v})^*P\right){\scriptstyle[(i-1)p_2+j,(i-1)p_2+j]}+\left(P^*{\bf X}_{k}^{v}({\bf X}_{k}^{v})^*P\right){\scriptstyle[(p_1-i-1)p_2+p_2-j,(p_1-i-1)p_2+p_2-j]}
\end{eqnarray*}
and
\begin{eqnarray*}
\frac{1}{\sqrt{p_1p_2}} \lambda{\scriptstyle[i,j]}({\bf X}_{k}^{v})\overline{\lambda{\scriptstyle[i,j]}({\bf X}_{k}^{v})}+\frac{1}{\sqrt{p_1p_2}} \lambda{\scriptstyle[p_1-i,p_2-j]}({\bf X}_{k}^{v})\overline{\lambda{\scriptstyle[p_1-i,p_2-j]}({\bf X}_{k}^{v})}
\end{eqnarray*}
are equal for any $1\leq i\leq p_1$ and $1\leq j\leq p_2$. Here, the entries of the matrix $\lambda(.)$ are taken modulo $p_1$ and $p_2$ and the entries of $[P^*{\bf X}_{k}^{v}({\bf X}_{k}^{v})^*P]$ are taken modulo $p_1p_2$. The result of Lemma \ref{lemme_transformee_fourier} follows.

\subsection{Proof of Proposition \ref{proposition_pente}}

\begin{proof}[Proof of Proposition \ref{proposition_pente}]
We only consider the anisotropic case, since  the proof for isotropic estimation is analogous.
For any model $m\in \mathcal{M}_1$, we define 
\begin{eqnarray*}
 \Delta(m,m') := \gamma_{n,p,p}\left(\widehat{\theta}_{m,\rho}\right)+\pen(m) - \gamma_{n,p,p}\left(\widehat{\theta}_{m',\rho}\right) - \pen(m') \ .
\end{eqnarray*}
We aim at showing that with large probability, the quantity $\Delta(m,m')$ is positive for all small dimensional models $m$. 
Hence, we would conclude that the dimension of $\widehat{m}$ is large.
In this regard, we bound the deviations of the differences 
\begin{eqnarray*}
\gamma_{n,p,p}\left(\widehat{\theta}_{m,\rho}\right) - \gamma_{n,p,p}\left(\widehat{\theta}_{m',\rho}\right)  & = &  \left[\gamma_{n,p,p}\left(\widehat{\theta}_{m,\rho}\right)-\gamma_{n,p,p}\left(\theta_{m,\rho}\right)\right]+ \bigg[\gamma_{n,p,p}\left(\theta_{m,\rho}\right) -\gamma_{n,p,p}\left(\theta \right)\bigg]\\  &+ &\left[\gamma_{n,p,p}(\theta)-\gamma_{n,p,p}\left(\widehat{\theta}_{m',\rho}\right)\right]\ .
\end{eqnarray*}

\begin{lemma}\label{concentration_premier_deuxieme}
Let $K_2$ be some universal constant that we shall define in the proof.
With probability larger than $3/4$, 
$$\gamma_{n,p,p}(\theta)-\gamma_{n,p,p}(\theta_{m,\rho})\leq \frac{K_2}{2}\rho^2\varphi_{\text{max}}(\Sigma)\frac{d_m\vee 1}{np^2}\,$$  and 
$$ \gamma_{n,p,p}\left(\theta_{m,\rho}\right)-\gamma_{n,p,p}\left(\widehat{\theta}_{m,\rho}\right) \leq\frac{K_2}{2}\rho^2\varphi_{\text{max}}(\Sigma)\frac{d_m}{np^2}\ $$
for all models $m\in\mathcal{M}_1$.
\end{lemma}

\begin{lemma}\label{concentration_troisieme}
Assume that $p$ is larger than some numerical constant $p_0$.
With probability larger than $3/4$, it  holds that
$$ \gamma_{n,p,p}(\theta)- \gamma_{n,p,p}(\widehat{\theta}_{m',\rho}) \geq K_3\sigma^2\left\{\varphi_{\text{min}}\left(I_{p^2}-C(\theta)\right)\wedge\left[\rho-\varphi_{\text{max}}\left(I_{p^2}-C(\theta)\right)\right] \right\}\frac{d_{m'}}{np^2}\  ,$$
where $K_3$ is a universal constant defined in the proof.
\end{lemma}
Let us take $K_1$ to be exactly $K_3$.  Gathering the two last lemma with  Assumption (\ref{hypothese_pente}), there exists an event $\Omega$ of probability larger than $1/2$ such that 
\begin{eqnarray*}
\Delta(m,m') \geq &\\ \frac{\sigma^2}{np^2} &\left\{K_1\eta d_{m'}\left[\varphi_{\text{min}}\left(I_{p^2}-C(\theta)\right)\wedge\left[\rho-\varphi_{\text{max}}\left(I_{p^2}-C(\theta)\right)\right] \right]-K_2\frac{(d_m\vee 1)\rho^2}{\varphi_{\text{min}}\left(I_{p^2}-C(\theta)\right)}\right\}\ , 
\end{eqnarray*}
for all models $m\in\mathcal{M}_1$.			
Thus, conditionally to $\Omega$, $\Delta(m,m')$ is positive  for all models $m\in\mathcal{M}_1$ that satisfy 
$$\frac{d_m\vee 1}{d_{m'}}\leq \frac{K_3\eta}{K_2\rho^2}\varphi_{\text{min}}\left(I_{p^2}-C(\theta)\right)\left\{\varphi_{\text{min}}\left(I_{p^2}-C(\theta)\right)\wedge\left[\rho-\varphi_{\text{max}}\left(I_{p^2}-C(\theta)\right)\right] \right\}\ .$$
By Lemma 8.7 in \cite{verzelen_gmrf_theorie}, the dimension $d_{m'}$ is larger than $0.5[\sqrt{np^2}\wedge (p^2-1)]$.
We conclude that $$d_{\widehat{m}_{\rho}}\vee 1\geq  \left[\sqrt{np^2}\wedge p^2-1\right]\frac{K_3\eta}{K_2\rho^2}\varphi_{\text{min}}\left(I_{p^2}-C(\theta)\right)\left\{\varphi_{\text{min}}\left(I_{p^2}-C(\theta)\right)\wedge\left[\rho-\varphi_{\text{max}}\left(I_{p^2}-C(\theta)\right)\right] \right\}\ ,$$
with probability larger than $1/2$.
\end{proof}\vspace{0.5cm}

\begin{proof}[Proof of Lemma \ref{concentration_premier_deuxieme}]
In the sequel, $\overline{\gamma}_{n,p,p}(.)$ denotes the difference $\gamma_{n,p,p}(.)$ -$\gamma(.)$. 
Given a model $m$, we consider the difference
\begin{eqnarray*}
\gamma_{n,p,p}\left(\theta\right) -\gamma_{n,p,p}\left(\theta_{m,\rho}\right)& = &\overline{\gamma}_{n,p,p}\left(\theta\right) -\overline{\gamma}_{n,p,p}\left(\theta_{m,\rho}\right) - l(\theta_{m,\rho},\theta)\ .
\end{eqnarray*}
Upper bounding the difference of $\gamma_{n,p,p}$ therefore  amounts to bounding the difference of  $\overline{\gamma}_{n,p,p}$. By definition of $\gamma_{n,p,p}$ and $\gamma$, it expresses as 
\begin{eqnarray*}
\overline{\gamma}_{n,p,p}\left(\theta\right) -\overline{\gamma}_{n,p,p}\left(\theta_{m,\rho}\right)
 & = & \frac{1}{p^2} tr\left\{\left[(I_{p^2}-C(\theta))^2
  -(I_{p^2}-C(\theta_{m,\rho}))^2\right]\left(\overline{\bf X^vX^{v*}}-\Sigma\right)\right\}\ .
\end{eqnarray*}
The matrices $\Sigma$, $(I_{p^2}-C(\theta))$, and $(I_{p^2}-C(\theta_{m,\rho}))$ are symmetric block circulant. By Lemma \ref{codiagonalisation}, they are jointly diagonalizable in the same orthogonal basis. If we note $P$ an orthogonal matrix associated to this basis, then $C(\theta_{m,\rho})$, $C(\theta)$, and $\Sigma$ respectively decompose in
\begin{eqnarray*}
 C(\theta_{m,\rho})=P^*D(\theta_{m,\rho})P\ ,\, C(\theta)=P^*D(\theta)P\ \, \text{and } \Sigma=P^*D(\Sigma)P\ ,
\end{eqnarray*}
where the matrices $D(\theta_{m,\rho})$, $D(\theta)$, and $D(\Sigma)$ are diagonal.
\begin{eqnarray}
\overline{\gamma}_{n,p,p}\left(\theta\right) & -&\overline{\gamma}_{n,p,p}\left(\theta_{m,\rho}\right)=\nonumber \\ &  & \frac{1}{p^2}tr\left\{\left(D(\theta_{m,\rho})-D(\theta)\right)\left[2I_{p^2}-D(\theta)-D(\theta_{m,\rho})\right]D_{\Sigma}\left(\overline{\bf YY^*}-I_{p^2}\right)\right\}\ ,\label{difference_gamma}
\end{eqnarray}
where the matrix ${\bf Y}$ is defined as $P\sqrt{\Sigma^{-1}}{\bf X}^vP^*$. Its components follow independent standard Gaussian distributions. Since the matrices involved in (\ref{difference_gamma}) are diagonal, Expression (\ref{difference_gamma}) is a linear combination of centered $\chi^2$ random variables. We  apply the following lemma to bound its deviations.
\begin{lemma}\label{concentration_chi}
Let $(Y_1,\ldots, Y_D)$ be i.i.d. standard Gaussian variables. Let $a_1,\ldots, a_D$ be fixed numbers. We set 
$$\|a\|_{\infty}:= \sup_{i=1,\ldots, D}|a_i|,\ \ \ \|a\|_2^2 := \sum_{i=1}^D a_i^2$$
Let $T$ be the random variable defined by
$$T :=  \sum_{i=1}^Da_i\left(Y_i^2-1\right)\ .$$
Then, the following deviation inequality holds for any positive $x$	
\begin{eqnarray*}
 \mathbb{P}\left[T\geq 2\|a\|_2\sqrt{x}+2\|a\|_{\infty}x\right]\leq e^{-x}\ .
\end{eqnarray*}
\end{lemma}
This result is very close to Lemma 1 of Laurent and Massart in \cite{laurent98}. The only difference lies in the fact that they constrain the coefficients $a_i$ to be non-negative. Nevertheless, their proof easily extends to our situation. Let us define the matrix $a$ of size $n\times p^2$ as
$$a^i[j]:= \frac{D_{\Sigma}{\scriptstyle[i,i]}\left(D(\theta_{m,\rho}){\scriptstyle[i,i]}-D(\theta){\scriptstyle[i,i]}\right)\left(2-D(\theta[i,i]-D(\theta_{m,\rho})[i,i]\right)}{np^2}\ ,$$
for any $1\leq i\leq n$ and any $1\leq j\leq p^2$. Since the matrices $I-C(\theta)$ and $I-C(\theta_{m,\rho})$ belong to the set $\Theta_{\rho}^+$, their largest eigenvalue is smaller than $\rho$. By Definition (\ref{definition_perte}) of the loss function $l(.,.)$,  $\|a\|_2\leq 2\rho\sqrt{\varphi_{\text{max}}( \Sigma) l(\theta_{m,\rho},\theta)/(np^2)}$ and $\|a\|_{\infty}\leq 4\rho^2\varphi_{\text{max}} (\Sigma )/(np^2)$. By Applying Lemma \ref{concentration_chi} to Expression (\ref{difference_gamma}),
we conclude that 
$$\mathbb{P}\left[\overline{\gamma}_{n,p,p}\left(\theta\right) -\overline{\gamma}_{n,p,p}\left(\theta_{m,\rho}\right)\geq l(\theta_{m,\rho},\theta)	+12\rho^2\frac{\varphi_{\text{max}}(\Sigma )}{np^2}  x\right]\leq e^{-x}\ , $$
for any $x>0$. Consequently, for any $K>0$, the difference of $\gamma_{n,p,p}(.)$ satisfies
\begin{eqnarray*}
 \gamma_{n,p,p}(\theta)- \gamma_{n,p,p}(\theta_{m,\rho}) \leq  \frac{K}{2}\rho^2\varphi_{\text{max}} (\Sigma)\frac{d_{m}\vee 1}{np^2}\ ,
\end{eqnarray*}
simultaneously for all models $m\in\mathcal{M}_1$ with probability larger than $1-\sum_{m\in\mathcal{M}_1\setminus\emptyset}e^{-K(d_m\vee 1)/24}$. If $K$ is chosen large enough, the previous upper bound holds on an event of probability larger than $7/8$. Let us call $K'_2$ such a value.\\

Let us now turn to the second part of the result. As previously, we decompose the difference of empirical contrasts
\begin{eqnarray*}
\gamma_{n,p,p}\left(\theta_{m,\rho}\right) - \gamma_{n,p,p}\left(\widehat{\theta}_{m,\rho}\right) = \overline{\gamma}_{n,p,p}\left(\theta_{m,\rho}\right) - \overline{\gamma}_{n,p,p}\left(\widehat{\theta}_{m,\rho}\right) - l\left(\widehat{\theta}_{m,\rho},\theta_{m,\rho}\right)
\end{eqnarray*}
Arguing as in the proof of Theorem 3.1 in \cite{verzelen_gmrf_theorie}, we obtain an upper bound analogous to Eq.(49) in \cite{verzelen_gmrf_theorie}
\begin{eqnarray*}
 \overline{\gamma}_{n,p,p}\left(\theta_{m,\rho}\right) - \overline{\gamma}_{n,p,p}\left(\widehat{\theta}_{m,\rho}\right) \leq l(\widehat{\theta}_{m,\rho},\theta_{m,\rho}) + \rho^2\bigg\{\sup_{R\in \mathcal{B}^{\mathcal{H}'}_{m^2,m^2}}\frac{1}{p^2}tr\left[RD_{\Sigma}\left(\overline{\bf Y Y^*}-I_{p^2}\right)\right]\bigg\}^2\ .
\end{eqnarray*}
The set $\mathcal{B}^{\mathcal{H}'}_{m^2,m^2}$ is defined in the proof of Lemma 8.2 in \cite{verzelen_gmrf_theorie}. Its precise definition is not really of interest in this proof. Coming back to the difference of $\gamma_{n,p,p}(.)$, we get 
\begin{eqnarray*}
\gamma_{n,p,p}\left(\theta_{m,\rho}\right) - \gamma_{n,p,p}\left(\widehat{\theta}_{m,\rho}\right)\leq \rho^2\left\{\sup_{R\in \mathcal{B}^{\mathcal{H}'}_{m^2,m^2	}}\frac{1}{p^2}tr\left[RD_{\Sigma}\left(\overline{\bf YY^*}-I_{p^2}\right)\right]\right\}^2\ .
\end{eqnarray*}
We consecutively apply Lemma 8.3 and 8.4  in \cite{verzelen_gmrf_theorie} to bound the deviation of this supremum. Hence, for any positive number $\alpha$,
\begin{eqnarray}\label{majoration_combinaison_gamma_n}
\gamma_{n,p,p}\left(\theta_{m,\rho}\right) - \gamma_{n,p,p}\left(\widehat{\theta}_{m,\rho}\right)\leq  L_1(1+\alpha/2)\rho^2\varphi_{\text{max}}( \Sigma)\frac{ d_{m}}{np^2}\ .
\end{eqnarray}
with probability larger than $1- \exp[-L_2\sqrt{d_{m}}(\frac{\alpha}{\sqrt{1+\alpha/2}}\wedge \frac{\alpha^2}{1+\alpha/2})]$. Thus, there exists some numerical constant $\alpha_0$ such that the upper bound (\ref{majoration_combinaison_gamma_n}) with $\alpha=\alpha_0$ holds simultaneously for all models $m\in\mathcal{M}_1\setminus\emptyset$ with probability larger than $7/8$. Choosing $K_2$ to be the supremum of $K'_2$ and  $2L_1(1+\alpha_0/2)$ allows to conclude.

\end{proof}\vspace{0.5cm}

\begin{proof}[Proof of Lemma \ref{concentration_troisieme}]
Thanks to the definition (\ref{definition_gamman_alternative}) of $\gamma_{n,p,p}(.)$ we obtain
\begin{eqnarray*}
 \gamma_{n,p,p}(\theta)- \gamma_{n,p,p}(\widehat{\theta}_{m',\rho}) & = & \frac{1}{p^2}\sup_{\theta'\in \Theta_{m',\rho}^+}tr\left[\left(C(\theta')-C(\theta)\right)\left(2I_{p^2}- C(\theta)-C(\theta')\right)\Sigma\overline{{\bf Z}{\bf Z^*}}\right]\ ,
\end{eqnarray*}
where the $p^2\times n$ matrix ${\bf Z}$ is defined by ${\bf Z}:= \sqrt{\Sigma}^{-1}{\bf X^v}$. We recall that the matrices $\Sigma$, $C(\theta)$ and $C(\theta')$ commute since they are jointly diagonalizable by Lemma \ref{codiagonalisation}. 
Let $(\Theta^+_{m',\rho}-\theta)$ be the  set $\Theta^+_{m',\rho}$ translated by $\theta$.
Since $C(\theta)+C(\theta')=C(\theta+\theta')$, we lower bound the difference of $\gamma_{n,p,p}(.)$ as follows
\begin{eqnarray*}
 \gamma_{n,p,p}(\theta)- \gamma_{n,p,p}(\widehat{\theta}_{m',\rho})& = & \frac{1}{p^2}\sup_{\theta'\in \left(\Theta_{m',\rho}^+-\theta\right)} 2\sigma^2tr\left[C(\theta')\overline{{\bf Z}{\bf Z^*}}\right]  -  tr\left[C(\theta')^2\Sigma\overline{{\bf Z}{\bf Z^*}}\right]\\
&\geq & \frac{\sigma^2}{p^2}\sup_{\theta'\in \left(\Theta_{m',\rho}^+-\theta\right)}\left\{ 2tr\left[C(\theta')\overline{{\bf Z}{\bf Z^*}}\right]  -\varphi_{\text{min}}^{-1}\left[I_{p^2}-C(\theta)\right]  tr\left[C(\theta')^2\overline{{\bf Z}{\bf Z^*}}\right]\right\}\ .
\end{eqnarray*}
 Let us consider $\Psi_{i_1,j_1},\ldots, \Psi_{i_{d_{m'}},j_{d_{m'}}}$ a basis of the space $\Theta_{m'}$ defined in Eq.(14) of \cite{verzelen_gmrf_theorie}. Let $\alpha$ be a positive number that we shall define later. We then introduce $\theta'$ as
\begin{eqnarray*}
 \theta' := \varphi_{\text{min}}\left[I_{p^2}-C(\theta)\right]\frac{\alpha}{p^2}\sum_{k=1}^{d_{m'}}tr\left[C\left(\Psi_{i_k,j_k}\right)\overline{{\bf Z}{\bf Z^*}}\right]\Psi_{i_k,j_k}\ .
\end{eqnarray*}
Since $\theta$ is assumed to belong to $\Theta_{m',\rho}^+$, the parameter $\theta'$ belongs to $(\Theta_{m',\rho}^+-\theta)$ if $$\varphi_{\text{max}}[C(\theta')]\leq \varphi_{\text{min}}\left(I_{p^2}-C(\theta)\right)\hspace{0.4cm}\text{ and  } \hspace{0.4cm}\varphi_{\text{min}}[C(\theta')]\geq -\rho+ \varphi_{\text{max}}\left(I_{p^2}-C(\theta)\right)\ .$$.
 The largest eigenvalue of $C(\theta')$ is smaller than $\|\theta'\|_1$ whereas its smallest eigenvalue is larger than $-\|\theta'\|_1$.
Let us upper bound the $l_1$ norm of $\theta'$:
\begin{eqnarray}
 \|\theta'\|_1&= &2\varphi_{\text{min}}\left[I_{p^2}-C(\theta)\right]\frac{\alpha}{p^2}\sum_{k=1}^{d_{m'}}\left|tr\left[C\left(\Psi_{i_k,j_k}\right)\overline{{\bf Z}{\bf Z^*}}\right]\right|\nonumber\\
&\leq& 2\sqrt{\frac{\alpha}{p^2}\varphi_{\text{min}}\left[I_{p^2}-C(\theta)\right]d_{m'}tr\left[C(\theta')\overline{{\bf Z}{\bf Z^*}}\right]}\ . \label{majoration_l1}
\end{eqnarray} Hence, $\theta'$ belongs to $(\Theta_{m',\rho}^+-\theta)$ if 
\begin{eqnarray}\label{condition_theta'}
\|\theta'\|_1\leq \varphi_{\text{min}}\left(I_{p^2}-C(\theta)\right)\wedge \left[\rho-\varphi_{\text{max}}\left(I_{p^2}-C(\theta)\right)\right]\ .
\end{eqnarray} 
Thus, we get the lower bound
\begin{eqnarray}
 \gamma_{n,p,p}(\theta)- \gamma_{n,p,p}(\widehat{\theta}_{m',\rho}) & \geq & \frac{\sigma^2}{p^2}\left\{2tr\left[C(\theta')\overline{{\bf Z}{\bf Z^*}}\right] - \varphi^{-1}_{\text{min}}\left[I_{p^2}-C(\theta)\right]tr\left[C(\theta')^2\overline{{\bf Z}{\bf Z^*}}\right]\right\}\ , \label{condition1_pente}
\end{eqnarray}
as soon as  Condition (\ref{condition_theta'}) is  satisfied.

Let us now bound the deviations of the two random variables involved in (\ref{majoration_l1}) and (\ref{condition1_pente}) by applying Markov's and Tchebychev's inequality.  For the sake of simplicity, we assume that $d_{m'}$ is smaller than $(p^2-2p)/2$. In such a case, all the nodes in $m'$ are different from their symmetric in $\Lambda$. We omit the proof for $d_{m'}$ larger than $(p^2-2p)/2$ because the approach is analogous but the computations are slightly more involved.
  Straightforwardly, we get $$\mathbb{E}\left[tr\left(C(\theta')\overline{{\bf ZZ^*}}\right)\right] = 4\alpha\varphi_{\text{min}}\left[I_{p^2}-C(\theta)\right] \frac{d_{m'}}{n}\ ,$$
since the neighborhood $m'$ only contains points $(i,j)$ whose symmetric $(-i,-j)$ is different. A cumbersome but pedestrian computation leads to the upper bound
\begin{eqnarray*}
 \var\left[tr\left(C(\theta')\overline{{\bf ZZ^*}}\right)\right]\leq L_1\alpha^2\varphi^2_{\text{min}}\left[I_{p^2}-C(\theta)\right]\frac{d_{m'}}{n^2} \ ,
\end{eqnarray*}
where $L_1$ is a numerical constant.
Similarly, we upper bound the expectation of $tr\left[C(\theta')^2\overline{{\bf Z Z^*}}\right]$
\begin{eqnarray*}
 \mathbb{E}\left[tr\left(C(\theta')^2\overline{{\bf Z Z^*}}\right)\right] \leq L_2\alpha^2\varphi^2_{\text{min}}\left[I_{p^2}-C(\theta)\right]\frac{d_{m'}}{n} \ .
\end{eqnarray*}
Let us respectively apply Tchebychev's inequality and Markov's inequality to the variables $tr\left[C(\theta')\overline{{\bf ZZ^*}}\right]$ and $tr\left[C(\theta')^2\overline{{\bf Z Z^*}}\right]$. Hence, there exists an event $\Omega$ of probability larger than $3/4$ such that
\begin{eqnarray*}
2tr\left[C(\theta')\overline{{\bf ZZ^*}}\right]-\varphi^{-1}_{\text{min}}\left[I_{p^2}-C(\theta)\right]tr\left[C(\theta')^2\overline{{\bf Z Z^*}}\right]\geq \hspace{4cm}\\ \hspace{4cm} \varphi_{\text{min}}\left[I_{p^2}-C(\theta)\right] \frac{d_{m'}}{n}\bigg\{8\alpha\bigg(1-\sqrt{\frac{L'_1}{d_{m'}}}\bigg)-\alpha^2L'_2\bigg\}
\end{eqnarray*}
and 	
$$tr\left[C(\theta')\overline{{\bf ZZ^*}}\right]\leq 4\alpha \varphi_{\text{min}}\left[I_{p^2}-C(\theta)\right] \frac{d_{m'}}{n} \bigg(1+\sqrt{\frac{L'_1}{d_{m'}}}\bigg) \ .$$ 
In the sequel, we assume that $p$ is larger than some universal constant $p_0$, which ensures the dimension $d_{m'}$ to be larger than $4L'_1$. Gathering (\ref{majoration_l1}) with the upper bound on $tr\left[C(\theta')\overline{{\bf ZZ^*}}\right]$ yields
$$\|\theta'\|_1\leq 2\sqrt{2}\alpha\varphi_{\text{min}}\left[I_{p^2}-C(\theta)\right] \frac{d_{m'}}{\sqrt{np^2}}\leq  2\sqrt{2}\alpha\varphi_{\text{min}}\left[I_{p^2}-C(\theta)\right] \ ,$$
since $d_{m'}\leq p\sqrt{n}$. If $2\sqrt{2}\alpha$ is smaller than $1\wedge\left\{\left[\rho-\varphi_{\text{max}}\left(I_{p^2}-C(\theta)\right)\right\} \varphi^{-1}_{\text{min}}\left[I_{p^2}-C(\theta)\right]\right]$, then Condition (\ref{condition_theta'}) is fulfilled on the event $\Omega$ and it follows from (\ref{condition1_pente}) that 
$$\mathbb{P}\left\{\gamma_{n,p,p}(\theta)- \gamma_{n,p,p}(\widehat{\theta}_{m',\rho})  \geq  4\sigma^2\varphi_{\text{min}}\left[I_{p^2}-C(\theta)\right]\frac{d_{m'}}{np^2}\left[\alpha-\alpha^2L'_2/4\right]\right\}\geq \frac{3}{4}\ .$$
 Choosing $\alpha=\frac{2}{L'_2}\wedge \frac{\sqrt{2}}{4}\wedge\sqrt{2}\frac{\rho-\varphi_{\text{max}}\left(I_{p^2}-C(\theta)\right)}{ 4\varphi_{\text{min}}\left(I_{p^2}-C(\theta)\right)}$, we get
$$\mathbb{P}\left\{\gamma_{n,p,p}(\theta)- \gamma_{n,p,p}(\widehat{\theta}_{m',\rho})  \geq  K_3\sigma^2\left\{\varphi_{\text{min}}\left[I_{p^2}-C(\theta)\right]\wedge\left[\rho-\varphi_{\text{max}}\left(I_{p^2}-C(\theta)\right)\right] \right\}\frac{d_{m'}}{np^2}\right\}\geq \frac{3}{4}\ ,$$
where $K_3$ is an universal constant.
\end{proof}

\section*{Acknowledgements}
I am grateful to Pascal Massart and Liliane Bel for many fruitful discussions. I also thank the referees and the associate editor for their suggestions that led to an improvement of the manuscript.

\addcontentsline{toc}{section}{References}

\bibliographystyle{alpha}

\bibliography{spatial}

\end{document}